\documentclass[11pt,reqno]{amsart}
\usepackage[varg]{txfonts}

\usepackage{anyfontsize}

\usepackage{amssymb, amsmath, amsthm, amsbsy, mathrsfs, mathtools, mathdots, bbm, stackrel, cancel, ulem, graphicx, relsize, dynkin-diagrams, upgreek} 

\usepackage[toc, page]{appendix}

\usepackage[all, cmtip]{xy}
\usepackage{tikz, tikz-cd}

\usepackage[mathscr]{euscript}
\usepackage[scr=boondoxo]{mathalfa}

\usepackage{setspace}
\usepackage[shortlabels]{enumitem}
\setlist[enumerate]{itemsep=1pt, topsep=2pt}
\setlist[itemize]{itemsep=1pt, topsep=2pt}

\usepackage{color, hyperref}
\hypersetup{colorlinks=true, linktoc=all, linkcolor=blue, citecolor=blue}

\DeclareMathOperator*{\sym}{Sym}

\DeclareMathOperator{\op}{op}

\DeclareMathOperator{\rep}{Rep}

\DeclareMathOperator{\en}{End}

\DeclareMathOperator{\id}{Id}
\makeatletter
\newcommand{\oset}[3][0ex]{%
  \mathrel{\mathop{#3}\limits^{
    \vbox to#1{\kern-2\ex@
    \hbox{$\scriptstyle#2$}\vss}}}}
\makeatother

\newcommand{\tl}[1]{\tilde{#1}}

\newcommand{\ep}{\epsilon}
\newcommand{\ve}{\varepsilon}

\newcommand{\opl}{\mathlarger{\mathop{\oplus}}}

\newcommand{\A}{\alpha}
\newcommand{\B}{\beta}
\newcommand{\G}{\gamma}
\newcommand{\D}{\delta}
\newcommand{\la}{\lambda}

\newcommand{\To}{\Rightarrow}

\newcommand{\se}{\subseteq}

\newcommand{\C}{\mathbb{C}}

\newcommand{\Q}{\mathbb{Q}}

\newcommand{\Z}{\mathbb{Z}}

\newcommand{\tb}[1]{\textbf{#1}}

\newcommand{\ti}[1]{\textit{#1}}
\newcommand{\ul}[1]{\underline{#1}}
\newcommand{\ol}[1]{\overline{#1}}

\newcommand{\mr}[1]{\mathrm{#1}}
\newcommand{\cl}[1]{\mathcal{#1}}
\newcommand{\fk}[1]{\mathfrak{#1}}

\newcommand{\p}[1]{\begin{proof}#1\end{proof}}

\newcommand{\eq}[1]{\begin{equation}\begin{aligned}#1\end{aligned}\end{equation}}
\numberwithin{equation}{section}

\newcommand{\comment}[1]{}

\newcommand{\bee}[1]{$$\begin{aligned}#1\end{aligned}$$}

\theoremstyle{plain}
\newtheorem{thm}{Theorem}[section]
\newtheorem{lem}[thm]{Lemma}
\newtheorem{prp}[thm]{Proposition}

\newtheorem{con}[thm]{Conjecture}

\theoremstyle{definition}
\newtheorem{dfn}[thm]{Definition}

\theoremstyle{remark}
\newtheorem{rmk}[thm]{Remark}

\newcommand{\kd}{\textcolor{orange}}
\author{}
\date{\today}
\title{Intertwiners of representations of twisted quantum affine algebras}

\begin{document}

\author{Keshav Dahiya and Evgeny Mukhin} 
\address{EM: Department of Mathematical Sciences,
Indiana University Indianapolis,
402 N. Blackford St., LD 270, 
Indianapolis, IN 46202, USA}
\email{emukhin@iu.edu} 

\address{KD: Department of Mathematical Sciences,
Indiana University Indianapolis,
402 N. Blackford St., LD 270, 
Indianapolis, IN 46202, USA}
\email{kkeshav@iu.edu} 

\begin{abstract}
    We use the $q$-characters to compute explicit expressions of the $R$-matrices for first fundamental representations of all types of twisted quantum affine algebras.

\medskip

\centerline{
  }

  \centerline{
  \textbf{\textit{AMS Classification numbers:\ }}{16T25, 17B38, 18M15 (primary), 17B37, 81R12.}}
  
\end{abstract}

\maketitle

\section{Introduction}
The $R$-matrices corresponding to quantum affine Lie algebras $U_q\tl{\fk{g}}$ are central objects of the theory of integrable systems. 

Physically, the entries of the $R$-matrix can be interpreted as weights of $XXZ$-type models. The quantum Yang-Baxter equation (QYBE) satisfied by the $R$-matrix is the origin of integrability of these models.

Mathematically, an $R$-matrix is an intertwiner of tensor products of two irreducible $U_q\tl{\fk{g}}$ modules in two different orders. The $R$-matrix is a rational function of a spectral parameter (or of a spectral shift of one of the factors). The zeroes and poles of the $R$-matrix correspond to the values of the spectral parameter when the tensor product stops being irreducible and as a result the products in different orders stop being isomorphic.

The explicit $R$-matrices corresponding to first fundamental  modules $\tl{L}_{1}$ have been computed in many cases a long time ago. For (untwisted) classical types the $R$-matrix is given in \cite{J86}. For G$_2$ the $R$-matrix was computed in \cite{O86} and \cite{K90}. For other exceptional types (with the  omission of E$_8$) it is obtained in \cite{M90},  \cite{M91}, see also \cite{DGZ94}. In these cases, $\tl{L}_{1}^{\otimes 2}$ is multiplicity free as a $U_q\fk{g}$-module.

For E$_8$ the $R$-matrix is described in \cite{ZJ20} and \cite{DM25}.  

For twisted quantum affine algebras $U_q\tl{\fk{g}}^\sigma$, the formulas for the $R$-matrix corresponding to the first fundamental  modules in types A$^{(2)}_{2r}$, A$^{(2)}_{2r-1}$ are given in \cite{B85}, \cite{J86}, for D$_{r+1}^{(2)}$ in \cite{B85}, \cite{KKMMNN92}, and for D$_4^{(3)}$ in \cite{KMOY06}. 
There is also a formula for E$_6^{(2)}$ $R$-matrix in terms of a non-standard restriction to the algebra of finite type C$_4$ (as opposed to the F$_4$ obtained by removing the affine $0$-th node), see \cite{GMW96}. There are also $R$-matrices for some other modules in types A$^{(2)}_{2r}$, A$^{(2)}_{2r-1}$ and  D$_{r+1}^{(2)}$, see  \cite{DGZ96}. 

In all cases (except \cite{B85}, \cite{ZJ20}) the main tool is the Jimbo equation which is deduced from the commutativity of the $R$-matrix with $E_0$ generator.

In \cite{DM25}, we developed an alternative method to compute the $R$-matrix using the theory of $q$-characters. The $q$-characters give full information about submodules and quotient modules of $\tl{L}_{1}(z_1)\otimes \tl{L}_{1}(z_2)$ which allows us to compute poles of the $R$-matrix and the values of $z_1/z_2$ when the $R$-matrix is well-defined but non-invertible. Together with simple general properties of the $R$-matrix, see Lemma \ref{lemR}, it determines the $R$-matrix almost uniquely in the case the poles of the $R$-matrix are simple. In this paper we show that this method allows to recover the $R$-matrices for the first fundamental  modules of all twisted affine quantum algebras. In particular, our formula in the case of E$_6^{(2)}$ is new.

The twisted quantum affine algebras are much less studied, and we have to prove a number of technical results to apply our machinery, see Theorems \ref{cyclic tensor product}, \ref{q char arg}, \ref{poles thm}. Our main sources on twisted quantum affine algebras are papers  \cite{CP98}, \cite{H10}, and \cite{Da14}.  

The main part of the answer are two $3\times 3$ and four $2\times 2$
matrices corresponding to multiplicity $3$ and multiplicity $2$ summands in the decomposition of $\tl L_{1}^{\otimes 2}$, see Theorems \ref{thm:R Dt2}, \ref{thm:R E6t2}, and \ref{thm:R D4t3}. Note that these matrices depend on the choice of some vectors and our choice differs from that in  \cite{KKMMNN92}, \cite{KMOY06} as we prefer to work with the symmetric coproduct and bases which are orthonormal with respect to Shapovalov forms.

\medskip 

The paper is organized as follows.  In Section \ref{twisted quantum affine algebras} we recall the basics of twisted quantum affine algebras and their $q$-characters. In Section \ref{sec:trivial multiplicities.} we discuss the multiplicity free cases of A$_{2r-1}^{(2)}$ $(r\geq 3$) and A$_{2r}^{(2)}$ $(r\geq 1$). In Section \ref{sec:non-trivial multiplicities.} we give details of D$^{(2)}_{r+1}$ $(r\ge 2)$, E$_6^{(2)}$, and D$_4^{(3)}$. 

\section{Preliminaries}
\label{twisted quantum affine algebras}

\subsection{Twisted quantum affine algebras}

We use the following general notations.

\begin{enumerate}

\item Let $\fk{g}$ be the simple simply-laced finite-dimensional Lie algebra of type A$_{2r}\ (r\ge 1)$, A$_{2r-1}\ (r\ge 3)$, D$_{r+1}\ (r\ge 2)$ or E$_6$. Let $\mr{I}$ be the set of nodes of the Dynkin diagram of $\fk{g}$ and $C=(C_{ij})_{i,j\in\mr{I}}$ be the Cartan matrix of $\fk{g}$. We choose the numbering on these Dynkin diagrams as follows: 
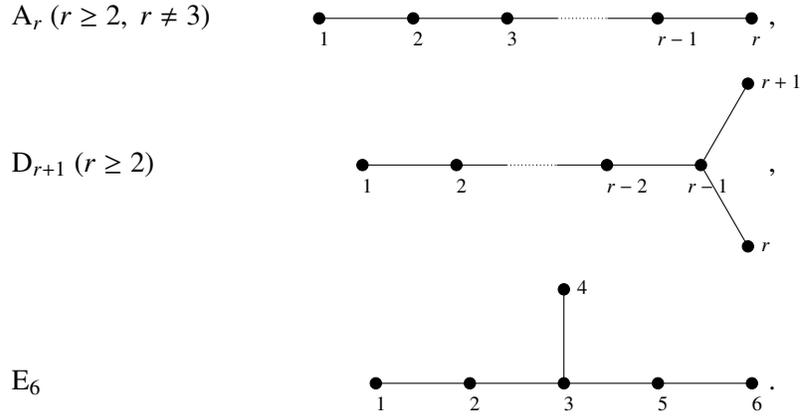
\bee{
& \text{A}_{r}\ (r\ge 2,\ r\ne 3) \qquad & \dynkin [edge length=1.25cm, root radius=0.075cm, label macro/.code={\drlap{#1}}, labels={1, 2, 3, r-1, r}, ordering=Kac, reverse arrows] A{***...**} \ ,\\
& \text{D}_{r+1}\ (r\ge 2) \qquad & \qquad \dynkin [edge length=1.25cm, root radius=0.075cm, label macro/.code={\drlap{#1}}, labels={1, 2, r-2, r-1, r, r+1}, ordering=Kac] D{} \hspace{-15pt} ,\\
& \text{E}_6\qquad \qquad & \dynkin [edge length=1.25cm, root radius=0.075cm, label macro/.code={\drlap{#1}}, labels={1, 2, 3, 5, 6, 4}, ordering=Kac] E{6} \ .\\
}

\item Let $\sigma$ be an automorphism of the Dynkin diagram of $\fk{g}$ of order $\mr{m}\in\{2,3\}$, that is, a bijection $\sigma:\mr{I}\to\mr{I}$ such that $\sigma\ne\id$, $\sigma^{\mr{m}}=\id$ and $C_{\sigma(i),\sigma(j)}=C_{i,j}$. Note that $\mr{m}=3$ only in the case of D$_4$. 
\item Let $\upomega\in\C$ be a primitive $\mr{m}$-th root of unity.

\item\label{identify Isigma} Let $\mr{I}^\sigma$ be the set of orbits of $\sigma$. For $i\in\mr{I}$ we denote the orbit of $i$ by $\ol{i}\in \mr{I}^\sigma$. We identify $\mr{I}^\sigma$ with
\begin{enumerate}
    \item $\{1,\dots,r\}$ in the case of A$_{2r}$, where $i$ is identified with the orbit $\{i,2r+1-i\}$.
    \item $\{1,\dots,r\}$ in the case of A$_{2r-1}$, where for $1\le i\le r-1$, $i$ is identified with the orbit $\{i,2r-i\}$, and $r$ is identified with the orbit $\{r\}$.
    \item $\{1,\dots,r\}$ in the case of D$_{r+1}$, where for $1\le i\le r-1$, $i$ is identified with the orbit $\{i\}$, and $r$ is identified with the orbit $\{r,r+1\}$.
    \item $\{1,2,3,4\}$ in the case of E$_6$, where for $1\le i\le 2$, $i$ is identified with the orbit $\{i,7-i\}$, $3$ is identified with the orbit $\{3\}$, and $4$ is identified with the orbit $\{4\}$.
    \item $\{1,2\}$ in the case of D$_4$ when $\mr{m}=3$, where $1$ is identified with the orbit $\{1,3,4\}$, and $2$ is identified with the orbit $\{2\}$.
\end{enumerate}
In all cases above we embed $\mr{I}^\sigma\subset \mr{I}$ as subsets of integers. This choice of embedding is fixed in what follows and we identify $\ol{i}\in \mr{I}^\sigma$ with $i\in \mr{I}$. For an $i\in\mr{I}^\sigma$, we say $i=\sigma(i)$ if the orbit of $i$ has cardinality one and $i\ne \sigma(i)$ if the orbit of $i$ has cardinality more than one.


\item The action of $\sigma$ on $\mr{I}$ is naturally extended to $\fk{g}$.  Let $\fk{g}^\sigma=\{g\in\fk{g}, \ \sigma(g)=g\}\subset\fk{g}$ be the Lie subalgebra fixed by $\sigma$. Then $\fk{g}^\sigma$ is  a simple finite-dimensional Lie algebra. Let  $C^\sigma=(C^\sigma_{ij})_{i,j\in\mr{I}^\sigma}$ be the corresponding Cartan matrix.


\item The action of $\sigma$ on $\fk{g}$ is extended to $\fk{g}\otimes \C[t,t^{-1}]$ by $\sigma(g\otimes t^k)=\sigma(g)\otimes (\upomega t)^k$, $g\in\fk{g}$. 
Let $\tl{\fk{g}}^\sigma=\{f\in\fk{g}\otimes \C[t,t^{-1}]:\sigma(f)=f\}$ be the Lie subalgebra fixed by $\sigma$. Then $\tl{\fk{g}}^\sigma$ is the twisted loop Lie algebra. Let $\tl{C}^\sigma=(\tl{C}^\sigma_{ij})_{i,j\in\tl{\mr{I}}^\sigma}$ be the corresponding Cartan matrix of the affine type. Here $\tl{\mr{I}}^\sigma=\{0\}\cup \mr{I}^\sigma$ and $\tilde C_{ij}^\sigma=C_{ij}^\sigma$ for $i,j\in \mr{I}^\sigma$.

\item The Cartan matrices $C^\sigma$ and $\tl{C}^\sigma$  can be read from the Dynkin diagrams given as follows: 
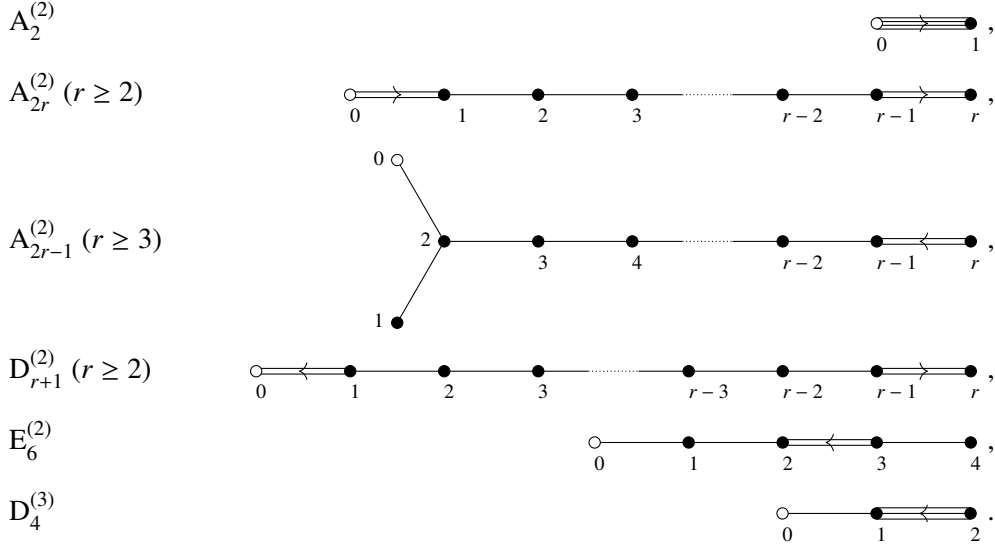
\bee{\allowdisplaybreaks
& \text{A}_{2}^{(2)}\qquad  & \dynkin [extended, edge length=1.25cm, root radius=0.075cm, label macro/.code={\drlap{#1}}, labels={0, 1}, ordering=Kac, reverse arrows] A[2]{2} \ ,\\
& \text{A}_{2r}^{(2)}\ (r\ge 2) \qquad & \dynkin [extended, edge length=1.25cm, root radius=0.075cm, label macro/.code={\drlap{#1}}, labels={0, 1, 2, 3, r-2, r-1, r}, ordering=Kac, reverse arrows] A[2]{even} \ ,\\
& \text{A}_{2r-1}^{(2)}\ (r\ge 3) \qquad & \dynkin [extended, edge length=1.25cm, root radius=0.075cm, label macro/.code={\drlap{#1}}, labels={0, 1, 2, 3, 4, r-2, r-1, r}, ordering=Kac] A[2]{odd} \ ,\\
& \text{D}_{r+1}^{(2)}\ (r\ge 2) \qquad & \dynkin [extended, edge length=1.25cm, root radius=0.075cm, label macro/.code={\drlap{#1}}, labels={0, 1, 2, 3, r-3, r-2, r-1, r}, ordering=Kac] D[2]{} \ ,\\
& \text{E}_6^{(2)}\qquad & \dynkin [extended, edge length=1.25cm, root radius=0.075cm, label macro/.code={\drlap{#1}}, labels={0, 1, 2, 3, 4}, ordering=Kac] E[2]{6} \ ,\\
& \text{D}_4^{(3)}\qquad & \dynkin [extended, edge length=1.25cm, root radius=0.075cm, label macro/.code={\drlap{#1}}, labels={0, 1, 2}, ordering=Kac] D[3]{4}\ .
}
\item Let $D^\sigma=\text{diag}\big(\{d_i\}_{i\in \mr{I}^\sigma}\big)$, respectively $\tl{D}^\sigma=\text{diag}\big(\{d_i\}_{i\in \tl{\mr{I}}^\sigma}\big)$, be such that $B^\sigma=D^\sigma C^\sigma$, respectively $\tl{B}^\sigma=\tl{D}^\sigma \tl{C}^\sigma$, is symmetric and $d_i\in\Z_{>0}$ are minimal possible except in the case of A$_{2r}^{(2)}$, $r\ge 1$, where $d_r=1/2$. The matrices $B^\sigma$, respectively $\tl{B}^\sigma$ are called the symmetrized Cartan matrices of finite, respectively affine, type.

\item Let $\A_i$, respectively $\omega_i$, $i\in \mr{I}^\sigma$, be simple roots, respectively fundamental weights, ${\cl{P}}=\oplus_{i\in\mr{I}^\sigma}\Z\omega_i$ the corresponding weight lattice and $\cl{P}_+=\oplus_{i\in\mr{I}^\sigma}\Z_{\ge0}\omega_i$ the cone of dominant weights.
We set $\omega_0=0\in\cl{P}_+$.

\item Let $a=(a_0,\dots,a_r)$ be the unique sequence of positive integers such that $\tl{C}^\sigma a^t=0$ and such that $a_0,\dots, a_r$ are relatively prime.

\item Let $q\in\C^\times$ be such that $q$ is not a root of unity. We fix a square root $q^{1/2}$.
Let $q_j=q^{d_j}$, $j\in \tl{\mr{I}}^\sigma$. For $k\in\frac{1}{2}\Z$ and $n\in\Z$, set
$$[n]_k=\frac{q^{kn}-q^{-kn}}{q^k-q^{-k}}\,\,,\quad [n]_k^{\mr{i}}=\frac{q^{kn}+(-1)^{n-1}q^{-kn}}{q^k+q^{-k}}.$$ 
Both $[n]_k$ and $[n]_k^{\mr{i}}$ are Laurent polynomials in $q^{1/2}$. We write $[n]_1$ as $[n]$ and $[n]_1^{\mr{i}}$ as $[n]^{\mr{i}}$.\\
Note that $\lim_{q\to 1}\,[n]_k=n$, $\lim_{q\to 1}[n]_k^{\mr{i}}=1$ if $n$ is odd, and $\lim_{q\to 1}[n]_k^{\mr{i}}=0$ if $n$ is even.

\item All representations are assumed to be finite-dimensional. We consider twisted quantum affine algebras of level zero only. All representations are assumed to be of type 1.
\end{enumerate}

\begin{dfn}[Drinfeld-Jimbo realization]\hspace{-10pt} \footnote{We follow \cite{H10}. In particular, our $U_q\tl{\fk{g}}^\sigma$ matches with the algebra in Definition 2 of \cite{Da14} for all types except of A$_{2r}^{(2)}$. In type A$_{2r}^{(2)}$ the algebra in \cite{Da14} coincides with our  $U_{q^2}\tl{\fk{g}}^\sigma$.} 
The twisted quantum affine algebra $U_q\tl{\fk{g}}^\sigma$ of level zero associated to $\fk{g}$ is an associative algebra over $\C$ with generators $E_i$, $F_i$, $K_i^{\pm1}$, $i\in\tl{\mr{I}}^\sigma$, and relations: 
$$K_iK_i^{-1}=K_i^{-1}K_i=1\,\,,\quad K_iK_j=K_jK_i\,\,,\quad K_0^{a_0}K_1^{a_1}\cdots K_r^{a_r}=1\,\,,$$
$$K_iE_jK_i^{-1}=q^{\tl{B}_{ij}^\sigma}E_j\,\,,\quad K_iF_jK_i^{-1}=q^{-\tl{B}_{ij}^\sigma}F_j\,\,,\quad [E_i,F_j]=\D_{ij}\frac{K_i-K_i^{-1}}{q_i-q_i^{-1}}\,\,,$$
$$\sum_{l=0}^{1-\tl{C}_{ij}^\sigma}(-1)^l\binom{1-\tl{C}_{ij}^\sigma}{l}_{q_i}E_i^lE_jE_i^{1-\tl{C}_{ij}^\sigma-l}=0\,\,,\quad\sum_{l=0}^{1-\tl{C}_{ij}^\sigma}(-1)^l\binom{1-\tl{C}_{ij}^\sigma}{l}_{q_i}F_i^lF_jF_i^{1-\tl{C}_{ij}^\sigma-l}=0\ ,\quad i\ne j\,\,.$$
\end{dfn}

The algebra $U_q\tl{\fk{g}}^\sigma$ has a Hopf algebra structure with comultiplication $\Delta$ 
\eq{\label{coproduct}
\Delta(K_i)=K_i\otimes K_i\,\,,\quad\Delta(E_i)=E_i\otimes K_i^{1/2}+K_i^{-1/2}\otimes E_i\,\,,\quad\Delta(F_i)=F_i\otimes K_i^{1/2}+K_i^{-1/2}\otimes F_i\,\,,\,\,\,\,i\in\tl{\mr{I}}^\sigma.
}
The antipode $\mr{S}$ given on the generators by 
\eq{\label{antipode}
\mr{S}(K_i)=K_i^{-1}\ ,\quad \mr{S}(E_i)=-K_i^{1/2} E_i K_i^{-1/2}\ ,\quad \mr{S}(F_i)=-K_i^{1/2}F_i K_i^{-1/2}\ ,\quad i\in\tl{\mr{I}}^\sigma\ .
}

The Hopf subalgebra of $U_q\tl{\fk{g}}^\sigma$ generated by $K_i^{\pm1}$, $E_i$, $F_i$, $i\in\mr{I}^\sigma$, is isomorphic to the quantum algebra $U_q\fk{g}^\sigma$ of finite type associated to $\fk{g}^\sigma$.

In what follows we also use the notation $U_q($B$_r)$,   $U_q($ F$_4)$, $U_q($A$_{2r}^{(2)})$, etc., for quantum algebras $U_q\fk{g}$ of type B$_r$, F$_4$, (twisted) quantum affine algebra $U_q\tl{\fk{g}}^\sigma$ of type  A$_{2r}^{(2)}$, etc.

The subalgebras $U_q\fk{g}^\sigma$ of $U_q\tl{\fk{g}}^\sigma$ in each case are as follows:
\bee{
\begin{array}{*{11}c}
    U_q(\text{A}_{2}^{(2)}) & & U_q(\text{A}_{2r}^{(2)}) & & U_q(\text{A}_{2r-1}^{(2)}) & & U_q(\text{D}_{r+1}^{(2)}) & & U_q(\text{E}_6^{(2)}) & & U_q(\text{D}_4^{(3)}) \\
    \cup & , & \cup & , & \cup & , & \cup & , & \cup & , & \cup \\
    U_{q^{1/2}}(\text{A}_1) & & U_{q^{1/2}}(\text{B}_r) & & U_q(\text{C}_r) & & U_q(\text{B}_r) & & U_q(\text{F}_4) & & U_q(\text{G}_2)
\end{array}\ .
}

\begin{dfn}
For $i\in\mr{I}^\sigma$, let $\tl{d}_i$ be $1$ in the case of A$_{2r}^{(2)}$ and $d_i$ otherwise. 

For $i,j\in\mr{I}$, let $d_{ij}\in\Z$ be $1$ in type A$_{2r}^{(2)}$, and in other types let $d_{ij}$ be given by
$$d_{ij}=\begin{cases}
    d_i & \text{if } C_{i,\sigma(i)}=0 \text{ and }\sigma(j)\ne j\ ,\text{ or if }\sigma(i)=i \ ,\\
    \mr{m} & \text{otherwise}\ .
\end{cases}$$

For $i,j\in\mr{I}$, let $P_{ij}^\pm(z_1,z_2)\in\Q(z_1,z_2)$ be given by
$$P_{ij}^\pm(z_1,z_2)=\begin{cases}
    1 & \text{if } C_{i,\sigma(i)}=0 \text{ and }\sigma(j)\ne j\ ,\text{ or if }\sigma(i)=i \ ,\\
    \dfrac{z_1^{\mr{m}}-q^{\pm 2 \mr{m}}z_2^{\mr{m}}}{z_1-q^{\pm 2}z_2} & \text{otherwise}\ .
\end{cases}$$
In type A$_{2r}^{(2)}$, we replace $q$ with $q^{1/2}$ in the above definition of $P_{ij}^{\pm}(z_1,z_2)$.
\end{dfn}

\begin{thm}[\cite{Dr87}\,\cite{Da14} Drinfeld's new realization]\label{Drinfeld realization}
The algebra $U_q\tl{\fk{g}}^\sigma$ is isomorphic to the algebra with generators $X_{i,n}^{\pm}$ $(i\in\mr{I}, n\in\Z)$, $K_i^{\pm1}$ $(i\in{\mr{I}})$, $\Phi_{i,\pm s}^{\pm}$ $(i\in\mr{I}, s\in\Z_{>0})$, and relations: 
$$X_{\sigma(i)}(z)=X_i(\upomega z)\ ,\quad \Phi_{\sigma(i)}(z)=\Phi_i(\upomega z)\ ,\quad K_{\sigma(i)}^{\pm1}=K_i^{\pm1}\ ,$$
$$K_iK_i^{-1}=K_i^{-1}K_i=1\,\,,\quad [\Phi_i^\pm(z_1),\Phi_j^\pm(z_2)]=[\Phi_i^\pm(z_1),\Phi_j^\mp(z_2)]=0\,\,,$$
$$\bigg(\prod_{k=1}^{\mr{m}}(q^{\pm C_{i,\sigma^k(j)}}\,\upomega^kz_1-z_2)\bigg)\,\Phi_i^\epsilon(z_1)X_j^\pm(z_2)=\bigg(\prod_{k=1}^{\mr{m}}(\,\upomega^kz_1-q^{\pm C_{i,\sigma^k(j)}}\,z_2)\bigg)\,X_j^\pm(z_2)\Phi_i^\epsilon(z_1)\,\,\text{ for }\epsilon=\pm\,,$$
$$\bigg(\prod_{k=1}^{\mr{m}}(q^{\pm C_{i,\sigma^k(j)}}\,\upomega^kz_1-z_2)\bigg)\,X_i^\pm(z_1)X_j^\pm(z_2)=\bigg(\prod_{k=1}^{\mr{m}}(\,\upomega^kz_1-q^{\pm C_{i,\sigma^k(j)}}\,z_2)\bigg)\,X_j^\pm(z_2)X_i^\pm(z_1)\ ,$$
$$[X_{i}^+(z_1),X_{j}^-(z_2)]=\frac{1}{\tl{d}_{\ol{i}}}\,\sum_{k=1}^{\mr{m}}\D_{i,\sigma^k(j)}\,\D\bigg(\frac{z_1}{\upomega^kz_2}\bigg)\,\frac{\Phi_{i}^+(z_1)-\Phi_{i}^-(z_1)}{q_{\ol{i}}-q_{\ol{i}}^{-1}}\ ,$$ 
$$\sym_{z_1,z_2}\bigg\{P_{ij}^\pm(z_1,z_2)\Big(X_j^\pm(z)X_i^\pm(z_1)X_i^\pm(z_2)-[2]_{d_{ij}}X_i^\pm(z_1)X_j^\pm(z)X_i^\pm(z_2)+X_i^\pm(z_1)X_i^\pm(z_2)X_j^\pm(z)\Big)\bigg\}=0,\ C_{ij}=-1,\sigma(i)\ne j\ ,$$
$$\sym\limits_{z_1,z_2,z_3}\bigg\{\Big(q^{3\ep/2}z_1^{\pm\ep}-[2]_{1/2}\,z_2^{\pm\ep}+q^{-3\ep/2}z_3^{\pm\ep}\Big)X_i^\pm(z_1)X_i^\pm(z_2)X_i^\pm(z_3)\bigg\}=0,\ C_{i,\sigma(i)}=-1, \ep=\pm1\ ,$$
where $\Phi_i^{\pm}(z)\in U_q\tl{\fk{g}}[\![z^{\pm1}]\!]$, $X_i^{\pm}(z)\in U_q\tl{\fk{g}}[\![z,z^{-1}]\!]$ for $i\in \mr{I}$, and $\D(z)\in\C[\![z,z^{-1}]\!]$ are given by
$$\Phi_i^\pm(z)=K_i^{\pm1}\Big(1+\hspace{-2pt}\sum_{ s\in\Z_{>0}}\Phi_{i,\pm s}^\pm\, z^{\pm s}\Big)
\ ,\quad X_i^\pm(z)=\sum_{n\in\Z}X_{i,n}^\pm\,z^{n}\ ,\quad \D(z)=\sum_{i\in\Z}z^i\ .$$
\qed
\end{thm}

Here the generators $X_{i,0}^+$, $X_{i,0}^-$, $K_i$, $i\in\mr{I}$ are respectively mapped to Drinfeld-Jimbo generators $E_{\ol{i}}$, $F_{\ol{i}}$, $K_{\ol{i}}$.

\begin{prp}[The shift of spectral parameter automorphism $\tau_a$]
For any $a\in\C^\times$, there is a Hopf algebra automorphism $\tau_a$ of $U_q\tl{\fk{g}}^\sigma$ defined by 
$$\tau_a(X_i^\pm(z))=X_{i}^\pm(az)\,\,,\quad\tau_a(\Phi_{i}^\pm(z))=\Phi_{i}^\pm(az)\,\,,\quad i\in{\mr{I}}^\sigma\,.$$ 
\qed
\end{prp}

Given a $U_q\tl{\fk{g}}^\sigma$-module $V$ and $a\in\C^\times$, we denote by $V(a)$ the pull-back of $V$ by $\tau_a$.

\begin{dfn}[Weight space]
Given a $U_q\fk{g}^\sigma$-module $V$ and $\lambda=\sum_{i\in \mr{I}^\sigma}\lambda_i\omega_i\in\cl{P}$, define the subspace $V_\la\subset V$ of weight $\lambda$ by
$$V_\lambda=\{v\in V:K_iv=q_i^{\lambda_i}v, \ i\in\mr{I}^\sigma\}.$$
If $V_\lambda\ne0$, $\lambda$ is called a weight of $V$. A nonzero vector $v\in V_\lambda$ is called a vector of weight $\lambda$.
\end{dfn}
For every representation $V$ of $U_q\fk{g}^\sigma$ we have $V=\oplus_{\lambda}V_{\lambda}$.

\begin{dfn}[$\ell$-weight]
Given a $U_q\tl{\fk{g}}^\sigma$-module $V$ and $\G=(\G_i^\pm(z))_{i\in \mr{I}^\sigma}$, $\G_i^\pm(z)\in\C[\![z^{\pm1}]\!]$, a sequence of formal power series in $z^{\pm1}$, define the subspace of generalized eigenvectors of $\ell$-weight $\G$ to be
$$V[\G]=\{v\in V: (\Phi_i^\pm(z)-\G_i^\pm(z))^{\dim(V)}\, v=0, \ i\in\mr{I}^\sigma\}.$$
If $V[\G]\ne0$, $\G$ is called an $\ell$-weight of $V$. Note that for any $\ell$-weight 
$\gamma_i^+(0)\gamma_i^-(\infty)=1$.

For every representation $V$ of $U_q\tl{\fk{g}}^\sigma$ we have $V=\oplus_{\G}V[\G]$ and for every $\lambda\in\cl{P}$, $V_\lambda=\oplus_{\G} (V_\lambda\cap V[{\G}])$.

A non-zero vector $v$ is a vector of $\ell$-weight $\G$ if 
$$(\Phi_i^\pm(z)-\G_i^\pm(z))\ v=0, \ i\in\mr{I}^\sigma.$$

\end{dfn}

\begin{dfn}[Highest $\ell$-weight representations]
A nonzero vector $v$ of $\ell$-weight $\G$ in some $U_q\tl{\fk{g}}^\sigma$-module $V$ is called an $\ell$-singular vector if $$X_i^+(z)\,v=0\,,\ i\in{\mr{I}}^\sigma\,.$$ 
A representation $V$ of $U_q\tl{\fk{g}}^\sigma$ is called a highest $\ell$-weight representation if $V=U_q\tl{\fk{g}}^\sigma\,v$ for some $\ell$-singular vector $v$. In such case $v$ is called the highest $\ell$-weight vector.
\end{dfn}

Let $\cl{U}$ be the set of all $\mr{I}^\sigma$-tuples $ p=(p_i)_{i\in\mr{I}^\sigma}$ of polynomials $p_i\in\C[z]$, with constant term $1$. 

\begin{thm}[\cite{CP98}\,\cite{H10}]\label{classification theorem}
\leavevmode
\begin{enumerate}
\item Every irreducible representation of $U_q\tl{\fk{g}}^\sigma$ is a highest $\ell$-weight representation. 
\item Let $V$ be an irreducible representation of $U_q\tl{\fk{g}}^\sigma$ of highest $\ell$-weight $\big(\G_i^\pm(z)\big)_{i\in\mr{I}^\sigma}$. Then there exists an $\mr{I}^\sigma$-tuple $p=(p_i)_{i\in\mr{I}^\sigma}\in\cl{U}$ such that $$\G_i^\pm(z)=\begin{cases} q^{\mr{m}\deg(p_i)}\dfrac{p_i(z^{\mr{m}}q^{-\mr{m}})}{p_i(z^{\mr{m}}q^{\mr{m}})}\in\C[\![z^{\pm1}]\!] & \text{if } i=\sigma(i) \ ,\vspace{5pt}\\
q^{\deg(p_i)}\dfrac{p_i(zq^{-1})}{p_i(zq)}\in\C[\![z^{\pm1}]\!] & \text{if } i\ne \sigma(i)\ .\end{cases}$$
\item Assigning to $V$ the $\mr{I}^\sigma$-tuple $p\in\cl{U}$ defines a bijection between $\cl{U}$ and the set of isomorphism classes of irreducible representations of $U_q\tl{\fk{g}}^\sigma$.
\end{enumerate}
\qed
\end{thm}

The polynomials $p_i(z)$ are called \ti{Drinfeld polynomials}. We denote the irreducible 
$U_q\tl{\fk{g}}^\sigma$-module corresponding to Drinfeld polynomials $p$ by $\tl{L}_p$. 



\begin{dfn}[Fundamental representations]
For each $i\in \mr{I}^\sigma\subset \mr{I}$, let $\tl{L}_{i}=\tl{L}_{p^{(i)}}$ be the irreducible $U_q\tl{\fk{g}}^\sigma$-module corresponding to the Drinfeld polynomials given by
$$p^{(i)}=(1-\D_{ij}z)_{j\in\mr{I}^\sigma}\ .
$$    
For $i\in\mr{I}^\sigma$, we call $\tl{L}_{i}(a)$ the $i$-th fundamental representation of $U_q\tl{\fk{g}}^\sigma$.
\end{dfn}

The category $\fk{Rep}\big(U_q\tl{\fk{g}}^\sigma\big)$ of representations of $U_q\tl{\fk{g}}^\sigma$ is an abelian monoidal category.  Denote by $\rep U_q\tl{\fk{g}}^\sigma$ the
 Grothendieck ring of $\fk{Rep}\big(U_q\tl{\fk{g}}^\sigma\big)$.

\begin{thm}\label{cyclic tensor product}
    Let $i,j\in\mr{I}^\sigma$. The tensor product $\tl{L}_i(a)\otimes\tl{L}_j(b)$ is cyclic on the tensor product of highest weight vectors if $a/b\ne \upomega^l q^{k}$ where $l\in\Z$ and $k\in\Z_{\ge 0}$.
\p{
For all types except for A$_{2r}^{(2)}$ the theorem is proved as Theorem 3 in \cite{C00}.

In type  A$_{2r}^{(2)}$, the affine subalgebra corresponding to node $r$ is isomorphic to A$_{2}^{(2)}$. Therefore for $r$-th reflection in the Weyl group of type B$_r$,
one cannot use the $U_q\hat{\fk{sl}}_2$-result (see Theorem 2 in \cite{C00}). Moreover, such a result may be difficult to prove in general. However, $\tl{L}_i(a)$ as an  A$_{2}^{(2)}$-module contains irreducible components isomorphic to trivial module and to three-dimensional fundamental module only. For such modules we check that the required products are cyclic by a direct computation.
Namely, it is easy to show that $\tl{L}_1(a)\otimes\tl{L}_1(b)$ for $a/b=q^{-2}$ and $a/b=-q^{-3}$
are cyclic from the product of highest $\ell$-weight vectors. In all other cases, when $a/b=\pm q^{-k}$, $k>0$, $k\ne 2,3$, the tensor product is irreducible.

After that we follow the proof of Theorem 3 in \cite{C00}.
}
\end{thm}

\begin{rmk}
In fact, in all types except for A$_{2r}^{(2)}$, one can prove the analog of Theorem \ref{cyclic tensor product}  for tensor product $\tl{L}_{i_1}(a_1)\otimes\dots\otimes \tl{L}_{i_n}(a_n)$ of arbitrary number of fundamental modules by the same method as in \cite{C00}. We do not need that result in this paper.
\end{rmk}

The category $\fk{Rep}\big(U_q\fk{g}^\sigma\big)$ of representations of $U_q\fk{g}^\sigma$ is an abelian monoidal semi-simple category. We denote the corresponding Grothendieck ring by $\rep U_q\fk{g}^\sigma$. Irreducible modules in $\fk{Rep}\big(U_q\fk{g}^\sigma\big)$ are parameterized by integral dominant weights. For $\lambda\in \cl{P}_+$, 
denote the corresponding irreducible $U_q\fk{g}^\sigma$-module by $L_\lambda$.

The module $L_\la$ has a unique (up to a scalar) symmetric bilinear form $(\ ,\ )$, called Shapovalov form, such that $E_i^*=F_i$, $i\in \mr{I}^{\sigma}$. The Shapovalov form is non-degenerate.

In what follows we will choose a weighted basis of $L_{\omega_1}$ such that $E_i^T=F_i$, $i\in \mr{I}^{\sigma}$, where $T$ stands for transposition. This basis is automatically orthonormal with respect to the Shapovalov form (with an appropriate normalization of the latter), see Lemma 2.9 in  \cite{DM25}. 

\subsection{Twisted \texorpdfstring{$q$}--characters}

For each $i\in\mr{I}^\sigma$, $a\in \C^\times$, let\hspace{-2pt}
\footnote{The variables $Y_{i,a}$ in this paper correspond to $Z_{i,a^\mr{m}}$ in \cite{H10} whenever $i=\sigma(i)$ and to $Z_{i,a}$ whenever $i\ne\sigma(i)$. Moreover, $Y_{i,a}$ in this paper denotes an $\mr{I}^\sigma$-tuple of rational functions while in \cite{H10}, $Y_{i,a}$ and $Z_{i,a}$ denote $\mr{I}$-tuples of rational functions.
With the relation \eqref{Z relation} these two conventions are equivalent. 
}
$Y_{i,a}=\big((Y_{i,a})_j\big)_{j\in\mr{I}^\sigma}$
be an $\mr{I}^\sigma$-tuple of rational functions given by $(Y_{i,a})_j(z)=1$ if $i\neq j$ and 
$$(Y_{i,a})_i(z)=\begin{cases}
    q^{\mr{m}}\dfrac{1-q^{-\mr{m}}z^{\mr{m}}a^{\mr
    m}}{1-q^{\mr{m}}z^{\mr{m}}a^{\mr
    m}} & \text{if }i = \sigma(i)\ , \vspace{5pt}\\
    q\dfrac{1-q^{-1}za}{1-qza} & \text{if }i\ne \sigma(i)\ .
    \end{cases}$$
The $\mr{I}^\sigma$-tuple $Y_{i,a}$ is the highest $\ell$-weight of $\tl{L}_{i}(a)$. Note that we have a relation
\eq{\label{Z relation}
Y_{i,a}=Y_{i,\upomega a} \qquad  \textrm{if}\ \sigma(i)=i\ .
}

Let $\cl{Y}$ be the abelian group of $\mr{I}^\sigma$-tuples of rational functions generated by $\{Y_{i,a}^{\pm1}\}_{i\in\mr{I}^\sigma,\,a\in\C^\times}$ modulo relation \eqref{Z relation} with component-wise multiplication. 

By Theorem 3.2 in \cite{H10} (or, alternatively, see Proposition \ref{weights change by A} below) the $\ell$-weights of representations of $U_q\tl{\fk{g}}^\sigma$ belong to $\cl{Y}$.

\begin{dfn}[Dominant $\ell$-weights]
An $\ell$-weight is called dominant if it is a monomial in variables\\ $\{Y_{i,a}\}_{i\in\mr{I}^\sigma,\,a\in\C^\times}$. The set of dominant $\ell$-weights will be denoted by $\cl{Y}_+$.
\end{dfn}

By Theorem \ref{classification theorem} the dominant monomials are in a bijective correspondence with highest $\ell$-weights of irreducible $U_q\tl{\fk{g}}^\sigma$-modules. 
In other words, the semi-group $\cl{Y}_+$ is naturally identified with $\cl{U}$.

\begin{dfn}[$q$-character]
The $q$-character of a $U_q\tl{\fk{g}}^\sigma$-module $V$ is the formal sum $$\chi_q(V)=\sum_{\G\in\cl{Y}}\dim(V[\G])\,\G\,\,\in\Z_{\geq 0}[\cl{Y}].$$
\end{dfn}

\begin{thm}[\cite{H10}] \label{qchar injective}
The $q$-character map $\chi_q:\rep U_q\tl{\fk{g}}^\sigma\to\Z_{\geq 0}[\cl{Y}],$ sending $V\mapsto \chi_q(V)$, is an injective ring homomorphism.
\qed
\end{thm}

A $U_q\tl{\fk{g}}^\sigma$-module $V$ is called special if $\chi_q(V)$ contains a unique dominant monomial.

\begin{dfn}[Simple $\ell$-roots]
For each $i\in\mr{I}^\sigma$ and $a\in\C^\times$, let $A_{i,a}\in\cl{Y}$ be given by
$$A_{i,a}=\begin{cases}
    Y_{i,aq}Y_{i,aq^{-1}}\bigg(\prod\limits_{j\sim i,j=\sigma(j)}Y_{j,a}^{-1}\bigg)\bigg(\prod\limits_{j\sim i,j\ne\sigma(j)}\Big(\prod\limits_{k=1}^{\mr{m}}Y_{j,\upomega^ka}^{-1}\Big)\bigg) & \text{if } C_{i,\sigma(i)}=2\ ,\vspace{4pt}\\
    Y_{i,aq}Y_{i,aq^{-1}}\bigg(\prod\limits_{j\sim i}Y_{j,a}^{-1}\bigg) & \text{if } C_{i,\sigma(i)}=0\ , \vspace{4pt}\\
    
    Y_{i,aq}Y_{i,aq^{-1}}Y_{i,-a}^{-1}\bigg(\prod\limits_{j\sim i}Y_{j,a}^{-1}\bigg) & \text{if } C_{i,\sigma(i)}=-1\ .
\end{cases}$$
Here all products are over $j\in I^\sigma\subset I$, and for $i,j\in\mr{I}^\sigma$, we write  $j\sim i$ if and only if $C^\sigma_{ij}<0$.

We call $A_{i,a}$ a simple $\ell$-root of color $i$.
\end{dfn}

Denote $Y_{1,a}$ by $\mr{1}_a$ , $Y_{2,a}$ by $\mr{2}_a$ and so on. For $m_+\in \cl{Y}_+$, let $p(m_+)\in \cl{U}$ be the corresponding set of Drinfeld polynomials. Denote $\tl{L}_{p(m_+)}$ by $\tl L_{m_+}$, and $\chi_q(\tl{L}_{p(m_+)})$ by $\chi_q(m_+)$.

The simple roots are given explicitly for each case as follows. \begin{enumerate}[align=left]
    \item[\hspace{5pt}A$_2^{(2)}$:] $A_{1,a}=\mr{1}_{aq}\mr{1}_{aq^{-1}}\mr{1}_{-a}^{-1}$.
    \vspace{4pt}
    \item[\hspace{5pt}A$_{2r}^{(2)}$:] $A_{1,a}=\mr{1}_{aq}\mr{1}_{aq^{-1}}\mr{2}_a^{-1},\>\ A_{i,a}=\mr{i}_{aq}\mr{i}_{aq^{-1}}\mr{(i-1)}_a^{-1}\mr{(i+1)}_a^{-1}$, $2\le i\le r-1,\>\ A_{r,a}=\mr{r}_{aq}\mr{r}_{aq^{-1}}\mr{r}_{-a}^{-1}\mr{(r-1)}_a^{-1}$. 
    \vspace{4pt}
    \item[\hspace{5pt}A$_{2r-1}^{(2)}$:] $A_{1,a}=\mr{1}_{aq}\mr{1}_{aq^{-1}}\mr{2}_a^{-1},\>\ A_{i,a}=\mr{i}_{aq}\mr{i}_{aq^{-1}}\mr{(i-1)}_a^{-1}\mr{(i+1)}_a^{-1}$, $2\le i\le r-1,\>\ A_{r,a}=\mr{r}_{aq}\mr{r}_{aq^{-1}}\mr{(r-1)}_a^{-1}\mr{(r-1)}_{-a}^{-1}$.
    \vspace{4pt}
    \item[\hspace{5pt}D$_{r+1}^{(2)}$:] $A_{1,a}=\mr{1}_{aq}\mr{1}_{aq^{-1}}\mr{2}_a^{-1},\>\ A_{i,a}=\mr{i}_{aq}\mr{i}_{aq^{-1}}\mr{(i-1)}_a^{-1}\mr{(i+1)}_a^{-1}$, $2\le i\le r-2$, \vspace{2pt}\\
    $A_{r-1,a}=\mr{(r-1)}_{aq}\mr{(r-1)}_{aq^{-1}}\mr{(r-2)}_{a}^{-1}\mr{r}_a^{-1}\mr{r}_{-a}^{-1},\>\ A_{r,a}=\mr{r}_{aq}\mr{r}_{aq^{-1}}\mr{(r-1)}_{a}^{-1}$.
    \vspace{4pt}
    \item[\hspace{5pt}E$_{6}^{(2)}$:] $A_{1,a}=\mr{1}_{aq}\mr{1}_{aq^{-1}}\mr{2}_a^{-1},\>\ A_{2,a}=\mr{2}_{aq}\mr{2}_{aq^{-1}}\mr{1}_a^{-1}\mr{3}_{a}^{-1},\>\ A_{3,a}=\mr{3}_{aq}\mr{3}_{aq^{-1}}\mr{2}_{a}^{-1}\mr{2}_{-a}^{-1}\mr{4}_{a}^{-1},\>\ A_{4,a}=\mr{4}_{aq}\mr{4}_{aq^{-1}}\mr{3}_a^{-1}$.
    \vspace{4pt}
    \item[\hspace{5pt}D$_{4}^{(3)}$:] $A_{1,a}=\mr{1}_{aq}\mr{1}_{aq^{-1}}\mr{2}_{a}^{-1},\>\ A_{2,a}=\mr{2}_{aq}\mr{2}_{aq^{-1}}\mr{1}_{a}^{-1}\mr{1}_{\mr{j}a}^{-1}\mr{1}_{\mr{j}^2a}^{-1}$, where $\mr{j}$ is a primitive cube root of unity.
\end{enumerate}

\comment{
\begin{lem}[Expansion Lemma]
    Let $V$ be a $U_q\tl{\fk{g}}^\sigma$-module. Let $v\in V$ be of $\ell$-weight $\G=(\G_j)_{j\in\mr{I}}\in \cl{Z}$ and $\ol{i}\in\mr{I}^\sigma$ with representative $i\in\mr{I}$ be such that $X_{i,s}^+\,v=0$ for all $s\in \Z$. If 
    $$\G_i=\begin{cases}
        q^{k\mr{m}}\dfrac{1-aq^{-2k\mr{m}}z}{1-az} & \text{if }i=\sigma(i) 
        \vspace{5pt}\\
        q^{k}\dfrac{1-aq^{-2k}z}{1-az} & \text{if }i\ne \sigma(i)
    \end{cases}\ ,$$
    where $k\in \Z_{>0}$, then for any $r\in\Z$
\bee{
X_{i,s}^+\,\big(X_{i,r}^-\,v-a^r v_a)=0\ ,\quad s\in \Z\ ,
}
for some vector $v_a$ with $\ell$-weight $\G A_{\ol{i},a}^{-1}$. 
\p{
\kd{????}
}
\end{lem}
}

\begin{prp}[\cite{MY14}]\label{weights change by A} 
Let $V$ be a $U_q\tl{\fk{g}}^\sigma$-module and $i\in\mr{I}^\sigma$. Suppose $\G$ and $\G'$ are $\ell$-weights of $V$. Then 
\bee{
V_{\G'}\cap\opl_{r\in\Z}X_{i,r}^{\pm}\big(V_{\G}\big)\ne 0\ \To\ \G'=\G\,A_{i,a}^{\pm}\quad \text{for some }\ a\in\C^\times.
}
\end{prp}
\p{
The proof is similar to Proposition 3.8 in \cite{MY14}, except that for A$_{2r}^{(2)}$ when $i=r$, a few details are slightly different. In this case, when $j=i=r$, equation (3.6) in \cite{MY14}, modifies to 
\bee{
\big((1-q^2uz)(1+q^{-1}uz)\G'_i(u)-(q^2-uz)(q^{-1}+uz)\G_i(u)\big)\la(z)=0\ ,
}
where $\la(z)$ is a formal Laurent series in $z$. Then we have $\la(z)\sum_{n=0}^\infty\,u^n(b_n+c_nz+d_nz^2)=0$. This gives a second order recurrence relation on the series coefficients of $\la(z)$. There exists a nonzero solution if and only if for some $a\in\C^\times$, we have $b_n+c_na+d_na^2=0$ for all $n\in\Z_{\ge0}$. Further, this solution is unique since for $j\ne i$, as in \cite{MY14}, we have a first order recurrence relation which gives a unique solution. Therefore, there is some $a\in\C^\times$ such that
\bee{
\G'_i(u)\big(\G_i(u)\big)^{-1} = \frac{(q^2-ua)(q^{-1}+ua)}{(1-q^2ua)(1+q^{-1}ua)}=\big(A_{i,a}(u)\big)_i\ .
}
The case of $X_{i,r}^-$ is similar.
}

\begin{thm}\label{q char arg} 
Let $V$ be an irreducible  $U_q\tl{\fk{g}}^\sigma$-module. Let $m$ be an $i$-dominant monomial in $\chi_q(V)$ of multiplicity one  for some $i\in\mr{I}$. Let $b\in\C^\times$ and $m_-=m A^{-1}_{i,b}$. Suppose 
\begin{enumerate}
    \item \label{q char arg 0} The power of $Y_{i,bq^{-1}}$ in $m$ is not greater than the power of $Y_{i,b q}$ in $m$.

    \item\label{q char arg 1} $m A_{i,c}\not\in\chi_q(V)$ for all $c\in\C^\times$. 


     \item\label{q char arg 2} $m_-A_{j,c}\not\in\chi_q(V)$ for all $j\in \mr{I}$, $c\in\C^\times$ unless $(j,c)=(i,b)$.

     \item\label{q char arg 3} The multiplicity of $m_-$ in $\chi_q(V)$ is not greater than one.
\end{enumerate}
Then  multiplicity of $m_-$ in $\chi_q(V)$ is zero,  $m_-\notin \chi_q(V)$.
\end{thm}
\p{
The proof is same as in the untwisted case. See Theorem 2.14 in \cite{DM25}.
We note that condition \eqref{q char arg 0} is equivalent to asserting that the $i$-th component of the $\mr{I}^\sigma$-tuple of rational functions corresponding to $m$ has no pole at $z=b^{-1}$.
}

We apply Theorem \ref{q char arg} to extract $\chi_q(V)$ from a known tensor product. In all our cases this tensor product has two dominant monomials, and we use Theorem \ref{q char arg} to prove that one of them is not in $\chi_q(V)$.  That enables us to identify $\chi_q(V)$.
Note that the conditions in Theorem \ref{q char arg} are written in combinatorial terms and therefore can be easily verified.

\subsection{\texorpdfstring{$R$}--matrices}

There is a quasitriangular structure on the Hopf algebra $U_q\tl{\fk{g}}^\sigma$, see \cite{KR90}, \cite{LS90}, \cite{Da98}.

\begin{prp}\label{R prop}
The Hopf algebra $U_q\tl{\fk{g}}^\sigma$ is almost cocommutative and quasitriangular, that is, there exists an invertible element $\fk{R}\in U_q\tl{\fk{g}}^\sigma\,\hat{\otimes}\,U_q\tl{\fk{g}}^\sigma$ of a completion of the tensor product, such that $$\Delta^{\op}(a)=\fk{R}\,\Delta(a)\,\fk{R}^{-1}\,,\,\,\,a\in U_q\tl{\fk{g}}\,,$$
where $\Delta^{\op}(a)=P\circ\Delta(a)$, $P$ is the flip operator, and
\eq{\label{quasitriangularity}
(\Delta\otimes\id)(\fk{R})=\fk{R}_{13}\fk{R}_{23}\,,\quad(\id\otimes\Delta)(\fk{R})=\fk{R}_{13}\fk{R}_{12}\,,\quad\fk{R}_{12}\fk{R}_{13}\fk{R}_{23}=\fk{R}_{23}\fk{R}_{13}\fk{R}_{12}\,.
}
\qed
\end{prp}

The element $\fk{R}$ is called the universal $R$-matrix of $U_q\tl{\fk{g}}^\sigma$.

The universal $R$-matrix has weight zero and homogeneous degree zero:
$$
(K_i\otimes K_i) \fk{R}=\fk{R} (K_i\otimes K_i), \qquad (\tau_z\otimes \tau_z) \fk{R}=\fk{R} (\tau_z\otimes \tau_z),\qquad i\in \tl{\mr{I}}^\sigma,\ z\in\C^\times.
$$

\begin{dfn}[The trigonometric $R$-matrix]
Let $V$ and $W$ be two representations of $U_q\tl{\fk{g}}^\sigma$ and $\pi_V$, $\pi_W$ be the respective representations maps. The map $$\tilde R^{V,W}(z)=(\pi_{V(z)}\otimes\pi_{W})(\fk{R}):V(z)\otimes W\to V(z)\otimes W$$ is called the $R$-matrix of $U_q\tl{\fk{g}}^\sigma$ evaluated in $V(z)\otimes W$. 
\end{dfn}

\begin{dfn}[Normalized $R$-Matrix]
Let $V$, $W$ be representations of $U_q\tl{\fk{g}}^\sigma$ with highest $\ell$-weight vectors $v$ and $w$ respectively. Denote by ${R}^{V,W}(z)\in\en(V\otimes W)$ the normalized $R$-matrix satisfying: 
$${R}^{V,W}(z)=f_{V,W}^{-1}(z)\,\tilde R^{V,W}(z)\,,$$ where $f_{V,W}(z)$ is the scalar function defined by $\tilde {R}^{V,W}(z)(v\otimes w)=f_{V,W}(z)\,v\otimes w.$ 
\end{dfn}

The map 
\eq{\label{R check}
\check{R}^{V,W}(z)=P\circ R^{V,W}(z):V(z)\otimes W\to W\otimes V(z)
}
(if it exists) is an intertwiner (or a homomorphism) of $U_q\tl{\fk{g}}^\sigma$-modules. If $V$, $W$ are irreducible, then the module $V(z)\otimes W$ is irreducible for all but finitely many $z\in\C^\times$. If for some $z$, the module $V(z)\otimes W$ is irreducible, then $W\otimes V(z)$ is also irreducible and the intertwiner is unique up to a constant. 

Equation \eqref{quasitriangularity} translates into the following lemma.
\begin{lem}\label{lem:qybe}
    Let $V_i$, $i=1,2,3$, be representations of $U_q\tl{\fk{g}}$. 
\begin{enumerate}
    \item \label{tqybe1}$\displaystyle R^{V_1,V_2}_{12}(z)\,R^{V_1,V_3}_{13}(zw)\,R^{V_2,V_3}_{23}(w)=R^{V_2,V_3}_{23}(w)\,R^{V_1,V_3}_{13}(zw)\,R^{V_1,V_2}_{12}(z)\,.$
    \item \label{tqybe} $\displaystyle\check{R}^{V_1,V_2}_{23}(z)\,\check{R}^{V_1,V_3}_{12}(zw)\,\check{R}^{V_2,V_3}_{23}(w)=\check{R}^{V_2,V_3}_{12}(w)\,\check{R}^{V_1,V_3}_{23}(zw)\,\check{R}^{V_1,V_2}_{12}(z)\,.$
\end{enumerate}
\qed
\end{lem}
The equations in Lemma \ref{lem:qybe} are called trigonometric QYBE.

The $R$-matrix $\check{R}^{V,W}(z)$ depends on the choice of the coproduct. In this paper we use coproduct $\Delta$ given by \eqref{coproduct}. Let $\mathfrak{R}_{\text{op}}$ be the universal $R$ matrix  corresponding to coproduct $\Delta^{\text{op}}$ and $\check{R}^{V,W}_{\text{op}}(z)$  be that $R$-matrix evaluated in $V(z)\otimes W$. Then $\mathfrak{R}_{\text{op}}=P\,\mathfrak{R}P$ and
\eq{\label{Rop}
\check{R}^{V,W}_{\text{op}}(z)=P(\pi_V\otimes \pi_W)\big((\tau_{z}\otimes 1)(\mathfrak{R}_{\text{op}})\big) = P\check{R}^{W,V}(z^{-1})P.
}
We collect a few general well-known properties of the $R$-matrices, cf. Lemma 2.19 in \cite{DM25}.
\begin{lem} \label{lemR}
Let $V_i$, $i=1,2$, be representations of $U_q\tl{\fk{g}}^\sigma
$. 
\begin{enumerate}
    \item  The normalized intertwiner $\check{R}^{V_1,V_2}(z)$ is a rational function of $z$. \label{try}
    \item If $V_1=\tl {L}_i(a)$ is fundamental, then $\check{R}^{V_1,V_1}(1)=\id$.
    \item $\check{R}^{V_1,V_2}(z;q)=P\check{R}^{V_2,V_1}(z^{-1};q^{-1})P$.
    \item\label{inversion relation} $\check{R}^{V_1,V_2}(z)\,\check{R}^{V_2,V_1}(z^{-1})=\id$.
    \item\label{R is self adjoint} $\check{R}^{V_1,V_2}(z)$ is self-adjoint with respect to the tensor Shapovalov form.
\end{enumerate}
\qed
\end{lem}

We use the following conjecture.

\begin{con}\label{conj:simple poles}
    Suppose $V(a)\otimes V$ has a single non-trivial submodule. Then the normalized $R$-matrix $\check{R}^{V,V}(z)$ has at most simple pole at $z=a$.
\end{con}
In general, one expects that the order of the pole at $z=a$ of a normalized $R$-matrix is at most one less than the number of irreducible subfactors. 

Note that in the trivial multiplicity case, we do not need that conjecture and instead we use the following lemma.

\begin{lem}\label{R0}
Let $V_1$, $V_2$ be irreducible representations of $U_q\tl{\fk{g}}^\sigma$ such that as $U_q\fk{g}^\sigma$-modules, $V_1$, $V_2$ are irreducible of highest weights $\la$, $\mu$ respectively. Suppose that the tensor product $L_\la\otimes L_\mu=\oplus_{\nu}\,L_\nu$ has trivial multiplicities. Then
\eq{\label{R check 0}
\check{R}^{V_1,V_2}(0)=\sum_{\nu}(-1)^{\nu}\,q^{(C(\nu)-C(\la+\mu))/2}P_\nu^{\la,\mu}\ ,
}
where $P_{\nu}^{\la,\mu}$ are projectors onto $L_{\nu}$, $(-1)^\nu=\pm 1$ is the eigenvalue of the flip operator $P$ on the $q\to 1$ limit of $L_\nu$, and $C(\nu)=(\nu,\nu+2\rho)$, with $\rho$ being the half sum of all positive roots, and $(\ ,\ )$ be the standard scalar product given on simple roots by $(\A_i,\A_j)=B^\sigma_{ij}$.
\p{
The proof is similar to the one for untwisted case. See \cite{DGZ94} and \cite{DGZ96}. A few more details are given in Lemma 2.20 in \cite{DM25}.
}
\end{lem}

In the untwisted  cases it is known that the submodules of tensor products of fundamental modules correspond to zeroes and poles of $R$-matrices, see Theorem 6.7 of \cite{FM1}. We prove the corresponding statement in the twisted case.

\begin{thm}\label{poles thm}
The tensor product $\tl{L}_{i}(a)\otimes \tl{L}_{j}(b)$ of fundamental representations of $U_q\tl{\fk{g}}^\sigma$, is reducible if and only if the normalized $R$-matrix $\check{R}^{V,W}(z)$  where $V=\tl{L}_{i}(1)$, $W=\tl{L}_{j}(1)$, has a pole or a nontrivial kernel at $z=a/b$. In that case, $a/b$ is necessarily equal to $\upomega^l q^k$, where $l,k$ are integers.
\p{The if statement is trivial. We prove the only if part.

Suppose $\tl{L}_{i}(a)\otimes \tl{L}_{j}(b)$ is reducible.  Then the dual module $(\tl{L}_{i}(a)\otimes \tl{L}_{j}(b))^*$ is also reducible. (Recall the antipode, \eqref{antipode}.)

We observe that the module dual to  $\tl{L}_{i}(a)$ is isomorphic to $\tl{L}_{i}(c_h a)$, where $c_h$ is given by
$$\begin{matrix}
{\textrm {A}}_{2r}^{(2)} & {\textrm {A}}_{2r-1}^{(2)} & {\textrm {D}}_{r+1}^{(2)} &  {\textrm {E}}_{6}^{(2)}&  {\textrm {D}}_{4}^{(3)}
\vspace{5pt}\\
-q^{2r+1} & -q^{2r} & q^{2r} & -q^{12} & q^6
\end{matrix}\ .$$
The constant $c_h$ can be computed as the shift of the lowest monomial in the $q$-character of $\tl{L}_{i}(a)$ with respect to the top monomial.

Then the opposite tensor product  $\tl{L}_{j}(b)\otimes \tl{L}_{i}(a)$ is a shift of the dual module $(\tl{L}_{i}(a)\otimes \tl{L}_{j}(b))^*$. 

One of the modules $\tl{L}_{i}(a)\otimes \tl{L}_{j}(b)$ and $\tl{L}_{j}(b)\otimes \tl{L}_{i}(a)$ is cyclic from the tensor product of highest weight vectors. Then the other one is not cyclic from the tensor product of highest weight vectors since it has structure similar to the dual module. It follows that  these two modules are not isomorphic. Thus if the $R$-matrix $\check{R}^{V,W}(a/b)$ is well-defined it gives a $U_q\tl{\fk{g}}^\sigma$-module homomorphism between these two modules and thus it has to be degenerate.
}
\end{thm}
We note that our proof differs from that of \cite{FM1}.


The following lemma is used for the computation of the $R$-matrix in the cases with nontrivial multiplicity.

Let $V$ be the first fundamental representation of $U_q\tl{\fk{g}}^\sigma$. Then we choose a basis $\{v_i\}_{i=1}^d$ of $V$ with the following properties. Let $\bar v_i=v_{\bar i}=v_{d+1-i}$ if weight of $v_i$ is not zero and $\bar v_i=v_i$ otherwise. Then we require
that the sum of weights of $v_i$ and $\bar v_i$ is zero and that 

\eq{\label{choice of basis}
E_j v_i=\sum_{r} a_{ir}^{(j)} v_r \quad \text{ if\ and\ only\ if }\quad  F_j \bar v_i=\sum_{r} a_{ir}^{(j)} \bar v_r\ ,\quad j\in\mr{I}^\sigma\ .
}

We construct such a basis for each type by a direct computation.  
In fact, the basis we choose is also orthonormal with respect to the Shapovalov form, and we have $E_j^T=F_j$, $j\in \mr{I}^\sigma$.

Let $t:V\to V$ be a linear map such that $v_i\mapsto \bar v_i$. Note that $t^2=\id$.

\begin{lem}\label{R and flip commute}
Let $V$ be the first fundamental representation of $U_q\tl{\fk{g}}^\sigma$. Then
\eq{\label{t and Rop}
\check{R}^{V,V}(z)=(t\otimes t)P\check{R}^{V,V}(z)P(t\otimes t)\ .
}
Here $P$ is the flip operator.
\p{
The proof is same as in the untwisted cases. See Lemma 2.22 in \cite{DM25}.
}
\end{lem}

\section{Cases of trivial multiplicity}\label{sec:trivial multiplicities.}

From now on, $\check{R}(z)$ denotes the intertwiner $\check{R}^{\tl{L}_1,\tl{L}_1}(z):\tl{L}_1(az)\otimes \tl{L}_1(a)\to\tl{L}_1(a)\otimes\tl{L}_1(az)$. When it is necessary to emphasize the dependence on $q$ we write $\check{R}(z\,;q)$ in place of $\check{R}(z)$. 

The following matrix $\check{R}(z)$ for untwisted type $\fk{sl}_{r+1}$ quantum affine algebra, given in \cite{J86}, will be used in matrix unit formulas for twisted $R$-matrices of type A.
\eq{\label{RqA}
\check{R}(z)=\sum_{i=1}^{r+1} E_{ii}\otimes E_{ii}\,+\,\frac{z(q-q^{-1})}{q-q^{-1}z}\,\sum_{i<j}E_{ii}\otimes E_{jj}\,+\,\frac{q-q^{-1}}{q-q^{-1}z}\,\sum_{i>j}E_{ii}\otimes E_{jj}\,+\,\frac{1-z}{q-q^{-1}z}\,\sum_{i\ne j}E_{ij}\otimes E_{ji}\ .
}
Here $E_{ij}$ are matrix units corresponding to a chosen basis $\{v_i\}$ in each case, that is, $E_{ij}(v_k)=\D_{jk}v_i$.

For a space $L$, we denote $\cl{S}^2(L), \Lambda^2(L)\subset L\otimes L$ the symmetric and skew-symmetric squares of $L$.


\subsection{Type A\texorpdfstring{$_{2r-1}^{(2)},\ r\ge 3$}{2}}

The $2r$-dimensional $U_q($A$_{2r-1}^{(2)})$-module $\tl{L}_1(a)$ restricted to $U_q($C$_r)$ is isomorphic to $L_{\omega_1}$. As $U_q($C$_r)$-modules, we have  
\eq{\label{tensorAodt}
\underbracket[0.1ex]{L_{\omega_1}}_{2r}\otimes \underbracket[0.1ex]{L_{\omega_1}}_{2r}\cong \underbracket[0.1ex]{L_{2\omega_1}}_{r(2r+1)}\oplus \underbracket[0.1ex]{L_{\omega_2}}_{(r-1)(2r+1)}\oplus \underbracket[0.1ex]{L_{\omega_0}}_{1}\ .
}
In the $q\to 1$ limit, $L_{2\omega_1}\mapsto \cl{S}^2(L_{\omega_1})$ and $L_{\omega_2}\oplus L_{\omega_0}\mapsto \Lambda^2(L_{\omega_1})$.

The $q$-character of $\tl{L}_{\mr{1}_a}$ has $2r$ terms and there are no weight zero terms:
\bee{
\chi_q(\mr{1}_a) =\ & \mr{1}_a+\ul{\mr{1}_{aq^2}^{-1}\mr{2}_{aq}}+\cdots+\mr{(r-2)}_{aq^{r-1}}^{-1}\mr{(r-1)}_{aq^{r-2}}+\mr{(r-1)}_{aq^r}^{-1}\mr{r}_{aq^{r-1}} \\
& + \mr{(r-1)}_{-aq^r}\mr{r}_{aq^{r+1}}^{-1} + \mr{(r-2)}_{-aq^{r+1}}\mr{(r-1)}_{-aq^{r+2}}^{-1} + \cdots + \mr{1}_{-aq^{2r-2}} \mr{2}_{-aq^{2r-1}}^{-1} + \ul{\mr{1}_{-aq^{2r}}^{-1}}\ .
}
In our formulas for $q$-characters we underline non-dominant monomials $m$ which can produce dominant monomials after multiplication by $\chi_q(\mr{1}_b)$. Note that in such a case the dominant monomial has the form $m\mr{1}_b$.

Using the $q$-characters we compute the poles of $\check{R}(z)$ and the corresponding kernels and cokernels. 

\begin{lem}\label{poles Ato}
The poles of the $R$-matrix $\check{R}(z)$, the corresponding submodules, and quotient modules are given by
\begin{center}\begin{tabular}{c@{\hspace{1cm}} c c} Poles & Submodules & Quotient modules  \\

$q^{2}$ & $\tl{L}_{\mr{1}_a\mr{1}_{aq^{-2}}}\cong L_{2\omega_1}$ & $\hspace{20pt} \tl{L}_{\mr{2}_{aq^{-1}}}\cong L_{\omega_2}\oplus L_{\omega_0}$ 
\vspace{2pt}\\

$-q^{2r}$ & $\hspace{22pt} \tl{L}_{\mr{1}_a\mr{1}_{-aq^{-2r}}}\cong L_{2\omega_1}\oplus L_{\omega_2}$ & $\hspace{6pt} \tl{L}_{\scriptscriptstyle{1}}\cong L_{\omega_0}$
\end{tabular}\ .\end{center}
\p{
From the $q$-character $\chi_q(\mr{1}_a)$ we see that the additional dominant monomials in the product $\chi_q(\mr{1}_{a})\chi_q(\mr{1}_b)$ occur only for $a/b=q^{\pm 2}$ and $a/b=-q^{\pm 2r}$. For all other cases there is a unique dominant monomial $1_a1_b$ and therefore $\tilde L_{1_a}\otimes \tilde L_{1_b}$ is irreducible. 

For $a/b=q^{\pm 2}$ and $a/b=-q^{\pm 2r}$, we have exactly two dominant monomials. For example, if $a/b=q^2$, the two monomials are $1_a1_{aq^{-2}}$ and $2_{aq^{-1}}$.

 We use Theorem \ref{q char arg}, to show that the dominant monomial which is not of the form $\mr{1}_a\mr{1}_b$, does not belong to $\chi_q(\mr{1}_a\mr{1}_b)$. It follows that $\tilde L_{1_a}\otimes \tilde L_{1_b}$ is reducible.
 
 Then $\tl{L}_{\mr{1}_a\mr{1}_b}$ is  either a submodule of $\tl{L}_1(a)\otimes\tl{L}_1(b)$, or a quotient module.
If $a/b=q^{-2}$ or $a/b=-q^{-2r}$, then by Theorem \ref{cyclic tensor product}, $\tl{L}_{\mr{1}_a\mr{1}_b}$ is cyclic from the tensor product of highest weight vectors, hence a quotient module. If $a/b=q^{2}$ or $a/b=-q^{2r}$, using the duality as in the proof of Theorem \ref{poles thm}, we conclude $\tl{L}_{\mr{1}_a\mr{1}_b}$ is a submodule.

Finally, by Theorem \ref{poles thm}, we conclude that $z=q^2$, $z=-q^{2r}$, are poles of $\check{R}(z)$ and at $z=q^{-2}$, $z=-q^{-2r}$, $\check{R}(z)$ is well-defined but has a nontrivial kernel coinciding with $\tl{L}_{\mr{2}_{aq}}$ and $\tl{L}_{\scriptscriptstyle{1}}$ respectively .
}
\end{lem}

We choose a basis $\{v_i:1\le i\le 2r\}$ for $L_{\omega_1}$ so that 
$F_iv_i=v_{i+1}$ and $F_iv_{\ol{i+1}}=v_{\ol{i}}$,  where $\ol{i}=2r+1-i$, and $i=1,\dots,r$.
In the chosen basis, $v_1\otimes v_1$ is a singular vector of weight $2\omega_1$, and $q\,v_1\otimes v_2-v_2\otimes v_1$ is a singular vector of weight $\omega_2$. We generate respectively the modules $L_{2\omega_1}$ and $L_{\omega_2}$ using these singular vectors. 

Let 
$\ve_i^q=(-q)^{r+1-i}$, $\ve_{\ol{i}}^q=-\ve_i^{q^{-1}}$, $1\le i\le r$. A singular vector $v_0\in L_{\omega_1}^{\otimes2}$ of weight $\omega_0$ is given by 
$$v_0=\sum_{i=1}^{2r}\ve_i^qv_i\otimes v_{\ol{i}}\ .$$

For $\la=2\omega_1,\omega_2,\omega_0$, let $P_\la^q$ be the projector onto the $U_q(\text{C}_r)$-module $L_\la$ in the decomposition \eqref{tensorAodt}, and $E_{ij}$ be matrix units corresponding to the chosen basis, that is, $E_{ij}(v_k)=\D_{jk}v_i$.

\begin{thm}\label{Rq At odd}
In terms of projectors, we have
\eq{\label{R proj Aodt}
\check{R}(z)=P_{2\omega_1}^q-q^{-2}\frac{1-q^{2}z}{1-q^{-2}z}P_{\omega_2}^q-q^{-2r-2}\frac{(1-q^{2}z)(1+q^{2r}z)}{(1-q^{-2}z)(1+q^{-2r}z)}P_{\omega_0}^q\ .
}
In terms of matrix units, we have 
\eq{\label{RqAodt}
\check{R}(z)=\big(\check{R}(z)\big)_{\fk{sl}_{2r}}-\frac{(q-q^{-1})(1-z)}{(q-q^{-1}z)(q^{r}+q^{-r}z)}Q(z)\ ,
}
where $\big(\check{R}(z)\big)_{\fk{sl}_{2r}}$ is the $\fk{sl}_{2r}$ trigonometric $R$-matrix in \eqref{RqA} and $Q(z)$ is given by
\bee{
Q(z)=z\sum_{i+j<2r+1}\frac{\ve_i^q\ve_j^q}{q^{r}}\,E_{ij}\otimes E_{\ol{i}\,\ol{j}}
-\sum_{i+j>2r+1}\frac{\ve_i^q\ve_j^q}{q^{-r}}\,E_{ij}\otimes E_{\ol{i}\,\ol{j}}\,+\,\frac{q^{r-1/2}-q^{-r+1/2}z}{q^{1/2}+q^{-1/2}}\,\sum_{i+j=2r+1}E_{ij}\otimes E_{\ol{i}\,\ol{j}}\ .
}
\p{
The poles of the $R$-matrix are known by Lemma \ref{poles Ato}. Using Lemma \ref{R0}, we conclude that these poles are simple. For example, in case of the summand $L_{\omega_0}$ in \eqref{tensorAodt}, the poles of the corresponding rational function $f_{\omega_0}(z)$ are at $z=q^2$, $z=q^{2r}$. Then using $f_k(z)f_k(z^{-1})=1$ and $f_k(1)=1$, it must be that 
$$f_{\omega_0}(z)=-q^{-2m_1-2m_2r}\frac{(1-q^2z)^{m_1}(1+q^{2r}z)^{m_2}}{(1-q^{-2}z)^{m_1}(1+q^{-2r}z)^{m_2}},\quad m_1,m_2\in\Z_{\ge1}\ .$$ 
By Lemma \ref{R0} we have $\check{R}(0)P_{\omega_0}^q=-q^{-2r-2}P_{\omega_0}^q$. This gives $m_1=m_2=1$.
}
\end{thm}

One can directly check that the $R$-matrix commutes with the action of $E_0$ and $F_0$. Namely,
\eq{\label{R and E0 commutation relation}
\check{R}(a/b)\,\Delta E_0(a,b)=\Delta E_0(b,a)\,\check{R}(a/b)\quad\text{and}\quad \check{R}(a/b)\,\Delta F_0(a,b)=\Delta F_0(b,a)\,\check{R}(a/b)\ ,
}
where  
$$K_0=q^{-1}E_{11}+q^{-1}E_{22}+\sum_{i=3}^{2r-2}E_{ii}+qE_{2r-1,2r-1}+qE_{2r,2r}\ ,\quad E_0(a)=a\big(E_{2r-1,1}+E_{2r,2}\big)\ ,$$
and $F_0(a)$ is the transpose of $a^{-2}E_0(a)$.


In the rational case, after substituting $z=q^{2u}$ in \eqref{RqAodt} and taking the limit $q\to 1$, we obtain 
\eq{\label{RuAodt}
\check{R}(u)=\frac{1}{1-u}\big(I-uP\big)\ ,
}
which is the untwisted type A$_{2r-1}^{(1)}$ rational $R$-matrix.

\subsection{Type A\texorpdfstring{$_{2}^{(2)}$}{2}}

The $3$-dimensional $U_q($A$_2^{(2)})$-module $\tl{L}_1(a)$ restricted to $U_{q^{1/2}}($A$_1)$ is isomorphic to $L_{2\omega_1}$. \\
As $U_{q^{1/2}}($A$_1)$-modules we have 
\eq{\label{tensorA2t}
\underbracket[0.1ex]{L_{2\omega_1}}_{3}\otimes \underbracket[0.1ex]{L_{2\omega_1}}_{3}\cong \underbracket[0.1ex]{L_{4\omega_1}}_{5}\oplus \underbracket[0.1ex]{L_{2\omega_1}}_{3}\oplus \underbracket[0.1ex]{L_{\omega_0}}_{1}\ .
}
In the $q\to 1$ limit, $L_{4\omega_1}\oplus L_{\omega_0}\mapsto \cl{S}^2(L_{\omega_1})$ and $L_{2\omega_1}\mapsto \Lambda^2(L_{\omega_1})$.

The $q$-character of $\tl{L}_{\mr{1}_a}$ has $3$ terms and there is $1$ weight zero term (shown in box):
$$\chi_q(\mr{1}_a)=\mr{1}_a+\boxed{\ul{\mr{1}_{aq^2}^{-1}\mr{1}_{-aq}}}+\ul{\mr{1}_{-aq^3}^{-1}}\ .$$
Using the $q$-characters, we compute the poles of $\check{R}(z)$ and the corresponding kernels and cokernels.
\begin{lem}\label{poles At2}
The poles of the $R$-matrix $\check{R}(z)$, the corresponding submodules and quotient modules are given by
\begin{center}\begin{tabular}{c@{\hspace{1cm}} c c}
Poles & Submodules &  Quotient modules  \\

$q^{2}$ & $\hspace{7pt}\tl{L}_{\mr{1}_a\mr{1}_{aq^{-2}}}\cong L_{4\omega_1}\oplus L_{\omega_0}$ & $\hspace{20pt} \tl{L}_{\mr{1}_{-aq^{-1}}}\cong L_{2\omega_1}$ 
\vspace{2pt}\\

$-q^{3}$ & $\hspace{8pt} \tl{L}_{\mr{1}_a\mr{1}_{-aq^{-3}}}\cong L_{4\omega_1}\oplus L_{2\omega_1}$ & $\hspace{34pt}\tl{L}_{\scriptscriptstyle{1}}\cong L_{\omega_0}$ \\
\end{tabular}\ .\end{center}
\p{
The proof is similar to the proof of Lemma \ref{poles Ato}.
}
\end{lem}

We choose a basis $\{v_i:1\le i\le 3\}$ for $L_{2\omega_1}$ so that $F_iv_i=\sqrt{[2]_{1/2}}v_{i+1}$, $i=1,2$.
In the chosen basis, $v_1\otimes v_1$ is a singular vector of weight $2\omega_1$ and $q\,v_1\otimes v_2-v_2\otimes v_1$ is a singular vector of weight $\omega_2$. 
We generate respectively the modules $L_{2\omega_1}$ and $L_{\omega_2}$ using these singular vectors. 

Let $\ve_1^q=q^{1/2}$, $\ve_2^q=-1$ and $\ve_{3}^q=q^{-1/2}$. A singular vector $v_0\in L_{2\omega_1}^{\otimes2}$ of weight $\omega_0$ is given by 
$$v_0=\sum_{i=1}^{3}\ve_i^qv_i\otimes v_{\ol{i}}\ ,$$ 
where $\ol{i}=4-i$.

For $\la=4\omega_1,2\omega_1,\omega_0$, let $P_\la^q$ be the projector onto the $U_{q^{1/2}}($A$_1)$-module $L_\la$ in the decomposition \eqref{tensorA2t}, and $E_{ij}$ be matrix units corresponding to the chosen basis, that is, $E_{ij}(v_k)=\D_{jk}v_i$.

\begin{thm}
In terms of projectors, we have
\eq{\label{R proj A2t}
\check{R}(z)=P_{4\omega_1}^q-q^{-2}\frac{1-q^{2}z}{1-q^{-2}z}P_{2\omega_1}^q+q^{-3}\frac{1+q^{3}z}{1+q^{-3}z}P_{\omega_0}^q\ .
}
In terms of matrix units, we have 
\eq{\label{RqA2t}
\check{R}(z)=\big(\check{R}(z)\big)_{\fk{sl}_3}+\frac{(q-q^{-1})(1-z)}{(q-q^{-1}z)(q^{3/2}+q^{-3/2}z)}Q(z)\ ,
}
where $\big(\check{R}(z)\big)_{\fk{sl}_3}$ is the $\fk{sl}_{3}$ trigonometric $R$-matrix in \eqref{RqA} and $Q(z)$ is given by 
\bee{
Q(z)=z\sum_{i+j<4}\frac{\ve_i^q\ve_j^q}{q^{3/2}}E_{ij}\otimes E_{\ol{i}\,\ol{j}}-\sum_{i+j>4}\frac{\ve_i^q\ve_j^q}{q^{-3/2}}E_{ij}\otimes E_{\ol{i}\,\ol{j}}-\frac{q-q^{-1}z}{q^{1/2}+q^{-1/2}}\sum_{i+j=4,i\ne 2}E_{ij}\otimes E_{\ol{i}\,\ol{j}}-\frac{q^2-q^{-2}z}{q^{1/2}+q^{-1/2}}E_{22}\otimes E_{22}\ .
}
\p{
The poles of the $R$-matrix are known by Lemma \ref{poles At2}. Using Lemma \ref{R0}, as in the proof of Theorem \ref{Rq At odd}, we conclude that these poles are simple.
}
\end{thm}

One can directly check that the $R$-matrix commutes with the action of $E_0$ and $F_0$, where  
$$K_0=q^{-2}E_{11}+E_{22}+q^2E_{33}\ ,\quad E_0(a)=a\,E_{31}\ ,$$
and $F_0(a)$ is the transpose of $a^{-2}E_0(a)$.


In the rational case, after substituting $z=q^{2u}$ in \eqref{RqA2t} and taking the limit $q\to 1$, we obtain 
\eq{\label{RuA2t}
\check{R}(u)=\frac{1}{1-u}\big(I-uP\big)\ ,
}
which is the untwisted type A$_{2}^{(1)}$ rational $R$-matrix.

\subsection{Type A\texorpdfstring{$_{2r}^{(2)}$, $r\ge2$}{2}}

The $(2r+1)$-dimensional $U_q($A$_{2r}^{(2)})$-module $\tl{L}_1(a)$ restricted to $U_{q^{1/2}}($B$_r)$ is isomorphic to $L_{\omega_1}$. For $r>2$, as $U_{q^{1/2}}($B$_r)$-modules, we have 
\eq{\label{tensorAevt}
\underbracket[0.1ex]{L_{\omega_1}}_{2r+1}\otimes \underbracket[0.1ex]{L_{\omega_1}}_{2r+1}\cong \underbracket[0.1ex]{L_{2\omega_1}}_{r(2r+3)}\oplus \underbracket[0.1ex]{L_{\omega_2}}_{\binom{2r+1}{2}}\oplus \underbracket[0.1ex]{L_{\omega_0}}_{1}\ .
}
In the $q\to 1$ limit, $L_{2\omega_1}\oplus L_{\omega_0}\mapsto \cl{S}^2(L_{\omega_1})$ and $L_{\omega_2}\mapsto \Lambda^2(L_{\omega_1})$. For $r=2$, $L_{\omega_2}$ has to be replaced with $L_{2\omega_2}$.

The $q$-character of $\tl{L}_{\mr{1}_a}$ has $2r+1$ terms and there is $1$ weight zero term (shown in box):
$$\chi_q(\mr{1}_a)=\mr{1}_a + \ul{\mr{1}_{aq^2}^{-1}\mr{2}_{aq}} + \cdots + \mr{(r-1)}_{aq^r}^{-1}\mr{r}_{aq^{r-1}} + \boxed{\mr{r}_{aq^{r+1}}^{-1}\mr{r}_{-aq^r}} + \mr{(r-1)}_{-aq^{r+1}}\mr{r}_{-aq^{r+2}}^{-1} + \cdots + \mr{1}_{-aq^{2r-1}}\mr{2}_{-aq^{2r}}^{-1} + \ul{\mr{1}_{-aq^{2r+1}}^{-1}}\ .$$
Using the $q$-characters we compute the poles of $\check{R}(z)$ and the corresponding kernels and cokernels.
\begin{lem}\label{poles Ate}
The poles of the $R$-matrix $\check{R}(z)$, the corresponding submodules and quotient modules are given by
\begin{center}\begin{tabular}{c@{\hspace{1cm}} c c}
Poles & Submodules & Quotient modules  \\

$q^{2}$ & $\hspace{30pt} \tl{L}_{\mr{1}_a\mr{1}_{aq^{-2}}}\cong L_{2\omega_1}\oplus L_{\omega_0}$ & $\hspace{10pt} \tl{L}_{\mr{2}_{aq^{-1}}}\cong L_{\omega_2}$ 
\vspace{2pt}\\

$-q^{2r-1}$ & $\hspace{18pt} \tl{L}_{\mr{1}_a\mr{1}_{-aq^{-2r-1}}}\cong L_{2\omega_1}\oplus L_{\omega_2}$ & $\hspace{24pt} \tl{L}_{\scriptscriptstyle{1}}\cong L_{\omega_0}$
\end{tabular}\ .\end{center}
\p{
The proof is similar to the proof of Lemma \ref{poles Ato}.
}
\end{lem}

We choose a basis $\{v_i:1\le i\le 2r+1\}$ for $L_{\omega_1}$ so that $F_iv_i=v_{i+1}$,  $F_iv_{\ol{i+1}}=v_{\ol{i}}$, $i=1,\dots, r-1$, $\ol{i}=2r+2-i$, and for $i=r$, $F_r.v_r=\sqrt{[2]_{1/2}}v_{r+1}$, $F_r.v_{\ol{r+1}}=\sqrt{[2]_{1/2}}v_{\ol{r}}$. 
In the chosen basis, $v_1\otimes v_1$ is a singular vector of weight $2\omega_1$, and $q\,v_1\otimes v_2-v_2\otimes v_1$ is a singular vector of weight $\omega_2$. 
We generate respectively the modules $L_{2\omega_1}$ and $L_{\omega_2}$ using these singular vectors. 

Let $\ve_i^q=(-1)^{r+1-i}q^{r-i+1/2}$, $\ve_{\ol{i}}^q=\ve_i^{q^{-1}}$, $1\le i\le r$, $\ve_{r+1}^q=1$. A singular vector $v_0\in L_{\omega_1}^{\otimes2}$ of weight $\omega_0$ is given by 
$$v_0=\sum_{i=1}^{2r+1}\ve_i^qv_i\otimes v_{\ol{i}}\ .$$ 

For $\lambda=2\omega_1,\omega_2\, (2\omega_2\text{ when }r=2),\omega_0$, let $P_\la^q$ be the projector onto the $U_{q^{1/2}}(\text{B}_r)$-module $L_\la$ in the decomposition \eqref{tensorAevt}, and $E_{ij}$ be matrix units corresponding to the chosen basis, that is, $E_{ij}(v_k)=\D_{jk}v_i$.

\begin{thm}
In terms of projectors, we have 
\eq{\label{R proj Aevt}
\check{R}(z)=P_{2\omega_1}^q-q^{-2}\frac{1-q^{2}z}{1-q^{-2}z}P_{\omega_2}^q+q^{-2r-1}\frac{1+q^{2r+1}z}{1+q^{-2r-1}z}P_{\omega_0}^q\ .
}
In terms of matrix units, we have
\eq{\label{RqAevt}
\check{R}(z)=\big(\check{R}(z)\big)_{\fk{sl}_{2r+1}}+\frac{(q-q^{-1})(1-z)}{(q-q^{-1}z)(q^{r+1/2}+q^{-r-1/2}z)}Q(z)\ ,
}
where $\big(\check{R}(z)\big)_{\fk{sl}_{2r+1}}$ is the $\fk{sl}_{2r+1}$ trigonometric $R$-matrix in \eqref{RqA} and $Q(z)$ is given by
\bee{
Q(z)=z\sum_{i+j<2r+2}\frac{\ve_i^q\ve_j^q}{q^{r+1/2}}E_{ij}\otimes E_{\ol{i}\,\ol{j}}-\sum_{i+j>2r+2}\frac{\ve_i^q\ve_j^q}{q^{-r-1/2}}E_{ij}\otimes E_{\ol{i}\,\ol{j}}-\frac{q^r-q^{-r}z}{q^{1/2}+q^{-1/2}}\sum_{\substack{i+j=2r+2 \\ i\ne r+1}}E_{ij}\otimes E_{\ol{i}\,\ol{j}} \\ 
-\frac{q^{r+1}-q^{-r-1}z}{q^{1/2}+q^{-1/2}}E_{r+1,r+1}\otimes E_{r+1,r+1}\ .
}
Here, in the case of $r=2$, $P_{\omega_2}^q$ is replaced by $P_{2\omega_2}^q$.
\p{
The poles of the $R$-matrix are known by Lemma \ref{poles Ate}. Using Lemma \ref{R0}, as in the proof of Theorem \ref{Rq At odd}, we conclude that these poles are simple.
}
\end{thm}

Note that for $r=1$, the $\check{R}(z)$ in \eqref{R proj Aevt}, \eqref{RqAevt} respectively, reduces to the $\check{R}(z)$ in \eqref{R proj A2t}, \eqref{RqA2t} for A$_2^{(2)}$ case.

One can directly check that the $R$-matrix commutes with the action of $E_0$ and $F_0$, where 
$$K_0=q^{-2}\,E_{11}+\sum_{i=2}^{2r}E_{ii}+q^2E_{2r+1,2r+1}\ ,\quad E_0(a)=a\,E_{2r+1,1}\ ,$$
and $F_0(a)$ is the transpose of $a^{-2}E_0(a)$.



In the rational case, after substituting $z=q^{2u}$ in \eqref{RqAevt} and taking the limit $q\to 1$, we obtain 
\eq{\label{RuAevt}
\check{R}(u)=\frac{1}{1-u}\big(I-uP\big)\ ,
}
which is the untwisted type A$_{2r}^{(1)}$ rational $R$-matrix.

\section{Cases of non-trivial multiplicities}\label{sec:non-trivial multiplicities.}
The main results in this section are Theorems \ref{thm:R Dt2} , \ref{thm:R E6t2}, \ref{thm:R D4t3}, which give $R$-matrices of the first fundamental representations in types D$_{r+1}^{(2)}$, E$_6^{(2)}$, D$_4^{(3)}$. The proofs of these theorems are quite similar but contain a few straightforward calculations and we prefer to give those proofs in detail.

\subsection{Type D\texorpdfstring{$_{r+1}^{(2)}$, $r\ge2$}{2}}

The $(2r+2)$-dimensional $U_q\big($D$_{r+1}^{(2)}\big)$-module $\tl{L}_1(a)$ restricted to $U_q($B$_r)$ is isomorphic to $L_{\omega_1}\oplus L_{\omega_0}$. For $r>2$, as $U_q($B$_r)$-modules, we have
\eq{\label{tensorDt2}
\big(\tl{L}_1(a)\big)^{\otimes 2}\cong\big(\underbracket[0.1ex]{L_{\omega_1}}_{2r+1}\oplus \underbracket[0.1ex]{L_{\omega_0}}_{1}\big)^{\otimes 2}\cong \underbracket[0.1ex]{L_{2\omega_1}}_{r(2r+3)}\oplus \underbracket[0.1ex]{L_{\omega_2}}_{r(2r+1)}\oplus\, 2 \underbracket[0.1ex]{L_{\omega_1}}_{2r+1}\oplus\, 2 \underbracket[0.1ex]{L_{\omega_0}}_{1}\ .
}
In the $q\to 1$ limit, $L_{2\omega_1}\oplus L_{\omega_0}\mapsto \cl{S}^2(L_{\omega_1})$ and $L_{\omega_2}\mapsto \Lambda^2(L_{\omega_1})$. For $r=2$, $L_{\omega_2}$ has to be replaced with $L_{2\omega_2}$.

For $r=2$, the $q$-character of $\tl{L}_{\mr{1}_{a}}$ has $6$ terms and there are $2$ weight zero terms (shown in box):
$$\chi_q(\mr{1}_{a})=\mr{1}_{a}+\ul{\mr{1}_{aq^2}^{-1}\mr{2}_{aq}\mr{2}_{-aq}} + \boxed{\mr{2}_{-aq^3}^{-1}\mr{2}_{aq}} + \boxed{\mr{2}_{aq^3}^{-1}\mr{2}_{-aq}} + \mr{1}_{aq^2}\mr{2}_{aq^3}^{-1}\mr{2}_{-aq^3}^{-1} + \ul{\mr{1}_{aq^4}^{-1}}\ .$$
Recall that for the case of D$_{r+1}^{(2)}$ we have $\sigma(i)=i$ and  $\mr{i}_{a}=\mr{i}_{-a}$ for $1\le i\le r-1$, see \eqref{Z relation}.

Using the $q$-characters, we compute the poles of $\check{R}(z)$ and the corresponding kernels and cokernels.
\begin{lem}
The poles of the $R$-matrix $\check{R}(z)$, the corresponding submodules and quotient modules are given by
\begin{center}\begin{tabular}{c@{\hspace{1cm}} c c}
Poles & Submodules & Quotient modules  \\

$\pm q^{2}$ & $\tl{L}_{\mr{1}_{a}\mr{1}_{aq^{-2}}}\cong L_{2\omega_1} \oplus L_{\omega_1} \oplus L_{\omega_0}$ & $\hspace{25pt} \tl{L}_{\mr{2}_{aq^{-1}}\mr{2}_{-aq^{-1}}}\cong L_{2\omega_2}\oplus L_{\omega_1} \oplus L_{\omega_0}$ 
\vspace{2pt}\\

$\pm q^{4}$ & $\hspace{40pt} \tl{L}_{\mr{1}_{a}\mr{1}_{aq^{-4}}}\cong L_{2\omega_1}\oplus L_{2\omega_2} \oplus 2\,L_{\omega_1} \oplus L_{\omega_0}$ & $\tl{L}_{\scriptscriptstyle{1}}\cong L_{\omega_0}$ \\
\end{tabular}\ .\end{center}
\p{
The proof is similar to the proof of Lemma \ref{poles Ato}.
}
\end{lem}

For $r>2$, the $q$-character of $\tl{L}_{\mr{1}_{a}}$ has $2r+2$ terms and there are $2$ weight zero terms (shown in box): 
\bee{
\chi_q(\mr{1}_{a})=\ & \mr{1}_{a} + \ul{\mr{1}_{aq^2}^{-1}\mr{2}_{aq}} + \cdots + \mr{(r-2)}_{aq^{r-1}}^{-1}\mr{(r-1)}_{aq^{r-2}} +\mr{(r-1)}_{aq^{r}}^{-1}\mr{r}_{aq^{r-1}}\mr{r}_{-aq^{r-1}} + \boxed{\mr{r}_{-aq^{r+1}}^{-1}\mr{r}_{aq^{r-1}}} \\ 
& + \boxed{\mr{r}_{aq^{r+1}}^{-1}\mr{r}_{-aq^{r-1}}} + \mr{(r-1)}_{aq^{r}} \mr{r}_{aq^{r+1}}^{-1} \mr{r}_{-aq^{r+1}}^{-1} + \mr{(r-2)}_{aq^{r+1}} \mr{(r-1)}_{aq^{r+2}}^{-1} + \cdots + \mr{1}_{aq^{2r-2}} \mr{2}_{aq^{2r-1}}^{-1} + \ul{\mr{1}_{aq^{2r}}^{-1}}\ .
}
Using the $q$-characters, we compute the poles of $\check{R}(z)$ and the corresponding kernels and cokernels.
\begin{lem}
The poles  of the $R$-matrix $\check{R}(z)$, the corresponding submodules and quotient modules are given by
\begin{center}\begin{tabular}{c@{\hspace{1cm}} c c}
Poles & Submodules & Quotient modules  \\

$\pm q^{2}$ & $\tl{L}_{\mr{1}_{a}\mr{1}_{aq^{-2}}}\cong L_{2\omega_1} \oplus L_{\omega_1} \oplus L_{\omega_0}$ & $\hspace{25pt} \tl{L}_{\mr{2}_{aq^{-1}}}\cong L_{\omega_2}\oplus L_{\omega_1} \oplus L_{\omega_0}$ 
\vspace{2pt}\\

$\pm q^{2r}$ & $\hspace{33pt} \tl{L}_{\mr{1}_{a}\mr{1}_{aq^{-2r}}}\cong L_{2\omega_1}\oplus L_{\omega_2} \oplus 2\,L_{\omega_1} \oplus L_{\omega_0}$ & $\hspace{-16pt}\tl{L}_{\scriptscriptstyle{1}}\cong L_{\omega_0}$ \\
\end{tabular}\ .\end{center}
\p{
The proof is similar to the proof of Lemma \ref{poles Ato}.
}
\end{lem}

We choose a basis $\{v_i:1\le i\le 2r+1\}\cup\{v_{2r+2}\}$ for $L_{\omega_1}\oplus L_{\omega_0}$ so that $F_iv_i=v_{i+1}$,  $F_iv_{\ol{i+1}}=v_{\ol{i}}$, $i=1,\dots, r-1$, $\ol{i}=2r+2-i$, and for $i=r$, $F_r.v_r=\sqrt{[2]}v_{r+1}$, $F_r.v_{r+1}=\sqrt{[2]}v_{\ol{r}}$. The vectors $v_{r+1}$ and $v_{2r+2}$ are of weight zero, and the vector $v_{2r+2}$ is annihilated by all $E_i$'s and $F_i$'s.

A singular vector of weight $2\omega_1$, respectively $\omega_2$, is chosen to be  $v_1\otimes v_1$, respectively $q\,v_1\otimes v_2-q^{-1}\,v_2\otimes v_1$. 
We choose the two singular vectors of weight $\omega_1$ to be $v_1\otimes v_{2r+2}\in L_{\omega_1}\otimes L_{\omega_0}$ and $v_{2r+2}\otimes v_1\in L_{\omega_0}\otimes L_{\omega_1}$ respectively. We choose the two singular vectors of weight $\omega_0$ to be respectively
$$w_1=\sum_{i=1}^{2r+1}\ve_i^qv_i\otimes v_{\ol{i}}\in L_{\omega_1}^{\otimes 2}\qquad \text{and}\qquad w_2=v_{2r+2}\otimes v_{2r+2}\in L_{\omega_0}^{\otimes 2}\ ,$$ 
where $\ve_i^q=(-1)^{i-1}q^{2r-2i+1}$, $\ve_{\ol{i}}^q=\ve_i^{q^{-1}}$, $1\le i\le r$, $\ve_{r+1}^q=(-1)^{r}$. 

For $\lambda=2\omega_1,\omega_2\, (2\omega_2\text{ when }r=2), \omega_1, \omega_0$, let $P_\lambda^q$ be the projector onto the $U_q($B$_r)$-module $L_\lambda$ in the decomposition \eqref{tensorDt2}.

\begin{thm}
\label{thm:R Dt2}
In terms of projectors, we have
\eq{\label{R proj Dt}
\check{R}(z)=P_{2\omega_1}^q-q^{-4}\frac{1-q^{4}z^2}{1-q^{-4}z^2}P_{\omega_2}^q + \frac{q^{-2}\,f_{\omega_1}(z)}{(1-q^{-4}\,z^2)}\otimes P_{\omega_1}^q+\frac{q^{-2r-2}\,f_{\omega_0}(z)}{(1-q^{-4}z^2)(1-q^{-4r}z^2)}\otimes P_{\omega_0}^q\ ,
}
where in the case of $r=2$, $P_{\omega_2}^q$ is replaced by $P_{2\omega_2}^q$, and the matrices $f_{\omega_1}(z)$, $f_{\omega_0}(z)$ are given by 
\bee{
f_{\omega_1}(z)=\begin{bmatrix}\B\,z & 1-z^2
\vspace{0.25cm} \\
1-z^2 & \B\,z \end{bmatrix}\ ,\quad 
f_{\omega_0}(z)=\begin{bmatrix}q^{-2r}+\A\,z^2+q^{2r}z^4 & \B\,z\,(1-z^2)  
\vspace{0.25cm} \\
\G\,z\,(1-z^2) & q^{2r}+\A\,z^2+q^{-2r}\,z^4 \end{bmatrix}\ .
}
Here $\A=[2]_{2r+2}-[2]_{2r}-[2]_{2r-2},\ $ $\B=[2]\,[2]^{\mr{i}}$ and $\G=[2]_{2r-1}\,[2]_{2r+1}^{\mr{i}}$.

\p{
In the expression of $\check{R}(z)$, the rational functions corresponding to the first two summands in \eqref{tensorDt2} are determined completely using the $q$-characters. Let $g_1(z)$ and $g_2(z)$ be the $2\times 2$ matrices corresponding to the last two summands $L_{\omega_1}$ and $L_{\omega_0}$ respectively.

The $2\times 2$ matrix $g_1(z)$ is determined completely as follows. Using Lemma \ref{R0},
\eq{\label{dtg10}
g_1(0)=\begin{bmatrix} 0 & q^{-2} \\ q^{-2} & 0
\end{bmatrix}\ ,\quad g_1(\infty)=\begin{bmatrix} 0 & q^{2} \\ q^{2} & 0
\end{bmatrix}\ .
}
From the $q$-characters we know the poles of $g_1(z)$ and by Conjecture \ref{conj:simple poles} we presume that the poles are simple. Combining this and \eqref{dtg10} with $g_1(1)$ being zero on off-diagonal entries and that $g_1(z)$ commutes with the flip operator acting on singular vectors, see Lemma \ref{R and flip commute}, we get
\bee{
g_1(z)=\frac{q^{-2}\,f_{\omega_1}(z)}{1-q^{-4}z^2}\quad\text{where}\quad f_{\omega_1}(z)=\begin{bmatrix}
a\,z & (1-z)(1+b\,z) \\
(1-z)(1+b\,z) & a\,z 
\end{bmatrix}\ ,
}
From $g_1(1)=\id$ we have $a=[2]\,[2]^{\mr{i}}$. From the inversion relation $g_1(z)g_1(z^{-1})=\id$, we get $b=1$.

The $2\times 2$ matrix $g_2(z)$ is determined (up to a sign) as follows. Using Lemma \ref{R0},
\eq{\label{g20}
g_2(0)=\begin{bmatrix}
q^{-4r-2} & 0\\ 0 & q^{-2}
\end{bmatrix}\ ,\quad g_2(\infty)=\begin{bmatrix}
q^{4r+2} & 0\\ 0 & q^2
\end{bmatrix}\ .
}
From the $q$-characters we know the poles of $g_2(z)$ and by Conjecture \ref{conj:simple poles} we presume that the poles are simple. Combining this and \eqref{g20} with $g_2(1)$ begin zero on off-diagonal entries we get
\bee{
g_2(z)=\frac{q^{-2r-2}\,f_{\omega_0}(z)}{(1-q^{-4}z^2)(1-q^{-4r}z^2)}\ ,
}
where 
\bee{
f_{\omega_0}(z)=\begin{bmatrix}
q^{-2r}+\A_1 z+\A z^2+\A_2 z^3+q^{2r}z^4 & z(1-z)(\beta_1+\beta_2 z) \\ z(1-z)(\G_1+\G_2 z) & q^{2r} +\A'_1 z+\A' z^2 + \A'_2 z^3 +q^{-2r} z^4
\end{bmatrix}\ .
}
From the inversion relation $g_2(z)g_2(z^{-1})=\id$, we get
\bee{
\B_1=\B_2\ ,\quad \G_1=\G_2\ ,\quad \A'=\A\ ,\quad \A'_1=\A_2\ ,\quad \A'_2=\A_2\ ,\quad \A_2=-q^{4r}\A_1\ ,
}
so that 
\bee{
f_{\omega_0}(z)=\begin{bmatrix}
q^{-2r}+\A_1 z+\A z^2-q^{4r}\A_1 z^3 + q^{2r}z^4 & \beta\, z(1-z^2) \\ \G\, z(1-z^2) & q^{2r} - q^{4r}\A_1 z + \A z^2 + \A_1 z^3 +q^{-2r} z^4
\end{bmatrix}\ .
}
Since $g_2(1)$ is $1$ on the diagonal entries we have
\eq{\label{g2 at 1}
\A_1(1-q^{4r})+\A+[2]_{2r}=[2]_2^{\mr{i}}\,[2]^{\mr{i}}_{2r}\ .
}
From $g_2(z)g_2(z^{-1})=\id$, now we get
\eq{\label{g2 inv a}
\A_1\big(\A+[2]_{2r}\big)=0\ ,
}
\eq{\label{g2 inv b}
\B\G+q^{4r}\A_1^2-\A[2]_{2r}=[2]_{4r}+[2]_4\ .
}
Now, using \eqref{g2 at 1} and \eqref{g2 inv a} we get two solutions for each of $\A$ and $\A_1$, 
\bee{
\text{either }\A=[2]_{2r+2}-[2]_{2r}-[2]_{2r-2}\ ,\ \A_1=0\ ,\quad \text{or}\quad \A=-[2]_{2r}\ ,\ \A_1=-q^{-2r}\,[2]_2\ .
}
Finally, from \eqref{g2 inv b}, we have
\bee{\text{either }\B\G=[2][2]^{\mr{i}}[2]_{2r-1}[2]_{2r+1}^{\mr{i}}\ ,\quad \text{or}\quad \B\G=0\ .
}
From the choice of singular vectors $w_1\in L_{\omega_1}^{\otimes2}$ and $w_2\in L_{\omega_0}^{\otimes2}$, we have 
$$\G=\frac{(w_1,w_1)}{(w_2,w_2)}\B=\frac{[2]_{2r-1}[2]_{2r+1}}{[2]}\B\ ,$$ 
so that 
\bee{
\text{either }\B=\pm [2][2]^{\mr{i}}\ ,\ \G=\pm [2]_{2r-1}[2]_{2r+1}^{\mr{i}}\quad \text{ or }\quad \B=\G=0\ .
}
The solution in the second case here is not correct and does not satisfy QYBE. To reject this extra solution and to fix the correct sign of $\B$ (or $\G$) in the first case, we use the $E_0$-action. Namely, we apply both sides of the commutation relation in \eqref{R and E0 commutation relation} to $v_1\otimes v_{2r+2}$ and compare coefficients of $v_{2r+2}\otimes v_{2r+2}$ on both sides.

One directly checks that the $R$-matrix commutes with the action of $E_0$ and $F_0$, where  
$$K_0=q^{-2}\,E_{11} + \sum_{i=2}^{2r}E_{ii} + q^2\,E_{2r+1,2r+1}\ ,\quad E_0(a)=a\,\sqrt{[2]}\, \big(E_{2r+2,1}+E_{2r+1,2r+2}\big)\ ,$$
and $F_0(a)$ is the transpose of $a^{-2}E_0(a)$.
}
\end{thm}



In the rational case, we recover the untwisted type D$_{r+1}^{(1)}$ rational $R$-matrix in Corollary 4.13 in \cite{DM25} as follows. Let $\check{R}(u)$ be the rational $R$-matrix obtained after substituting $z=q^{2u}$ in \eqref{R proj Dt} and taking the $q\to 1$ limit. Let $T:\C^{2r+2}\to\C^{2r+2}$ be a linear map given by $T(v_i)=v_i$, for $1\le i\le r$, $T(v_i)=v_{i+1}$ for $r+2\le i\le 2r+1$, $T(v_{r+1})=v_{r+1}+\frac{1}{2}v_{r+2}$, and $T(v_{2r+2})=v_{r+1}-\frac{1}{2}v_{r+2}$ when $r$ is even, while $T(v_{2r+2})=\mr{i}v_{r+1}-\frac{\mr{i}}{2}v_{r+2}$, when $r$ is odd. Here $\mr{i}$ is the primitive second root of unity. Then the matrix $(T\otimes T)\check{R}(u)(T\otimes T)^{-1}$ is the untwisted type D$_{r+1}^{(1)}$ rational $R$-matrix.

\subsection{Type E\texorpdfstring{$_6^{(2)}$}{2}}

The $27$-dimensional $U_q\big($E$_6^{(2)}\big)$-module $\tl{L}_1(a)$ restricted to $U_q($F$_4)$ is isomorphic to $L_{\omega_1}\oplus L_{\omega_0}$.\\ 
As $U_q($F$_4)$-modules, we have
\eq{\label{tensorE6t}
\big(\underbracket[0.1ex]{L_{\omega_1}}_{26}\oplus \underbracket[0.1ex]{L_{\omega_0}}_{1}\big)^{\otimes 2}\cong \underbracket[0.1ex]{L_{2\omega_1}}_{324}\oplus \underbracket[0.1ex]{L_{\omega_2}}_{273}\oplus \underbracket[0.1ex]{L_{\omega_4}}_{52}\oplus\,3 \underbracket[0.1ex]{L_{\omega_1}}_{26}\oplus\,2 \underbracket[0.1ex]{L_{\omega_0}}_{1}\ .
}
In the $q\to1$ limit, $L_{2\omega_1}\oplus L_{\omega_1}\oplus L_{\omega_0}\mapsto \cl{S}^2(L_{\omega_1})$ and $L_{\omega_2}\oplus L_{\omega_4}\mapsto \Lambda^2(L_{\omega_1})$.

The $q$-character of $\tl{L}_{\mr{1}_a}$ has $27$ terms and there are $3$ weight zero terms (shown in box):
\bee{
\chi_q(\mr{1}_a)=\ & \mr{1}_a + \mr{1}_{a q^2}^{-1} \mr{2}_{a q} + \mr{2}_{a q^3}^{-1} \mr{3}_{a q^2} + \mr{2}_{-a q^3} \mr{3}_{a q^4}^{-1} \mr{4}_{a q^3} + \mr{2}_{-a q^3} \mr{4}_{a q^5}^{-1} + \mr{1}_{-a q^4} \mr{2}_{-a q^5}^{-1} \mr{4}_{a q^3} \\ 
& + \ul{\mr{1}_{-a q^6}^{-1}\mr{4}_{a q^3}} + \mr{1}_{-a q^4} \mr{2}_{-a q^5}^{-1} \mr{3}_{a q^4} \mr{4}_{a q^5}^{-1} + \mr{1}_{-a q^6}^{-1} \mr{3}_{a q^4} \mr{4}_{a q^5}^{-1} + \mr{1}_{-a q^4} \mr{2}_{a q^5} \mr{3}_{a q^6}^{-1} + \mr{1}_{-a q^6}^{-1} \mr{2}_{a q^5} \mr{2}_{-a q^5} \mr{3}_{a q^6}^{-1} \\ 
& + \mr{1}_{a q^6} \mr{1}_{-a q^4} \mr{2}_{a q^7}^{-1} + \boxed{\mr{1}_{-a q^6}^{-1} \mr{1}_{a q^6} \mr{2}_{a q^7}^{-1} \mr{2}_{-a q^5}} + \boxed{\ul{\mr{1}_{a q^8}^{-1} \mr{1}_{-a q^4}}} + \boxed{\mr{2}_{-a q^7}^{-1} \mr{2}_{a q^5}} + \mr{1}_{a q^8}^{-1} \mr{1}_{-a q^6}^{-1} \mr{2}_{-a q^5} \\
& + \mr{1}_{a q^6} \mr{2}_{a q^7}^{-1} \mr{2}_{-a q^7}^{-1} \mr{3}_{a q^6} + \mr{1}_{a q^8}^{-1} \mr{2}_{-a q^7}^{-1} \mr{3}_{a q^6} + \mr{1}_{a q^6} \mr{3}_{a q^8}^{-1} \mr{4}_{a q^7} + \mr{1}_{a q^8}^{-1} \mr{2}_{a q^7} \mr{3}_{a q^8}^{-1} \mr{4}_{a q^7} + \mr{1}_{a q^6} \mr{4}_{a q^9}^{-1} \\ 
& + \mr{1}_{a q^8}^{-1} \mr{2}_{a q^7} \mr{4}_{a q^9}^{-1} + \mr{2}_{a q^9}^{-1} \mr{4}_{a q^7} + \mr{2}_{a q^9}^{-1} \mr{3}_{a q^8} \mr{4}_{a q^9}^{-1} + \mr{2}_{-a q^9} \mr{3}_{a q^{10}}^{-1} + \mr{1}_{-a q^{10}} \mr{2}_{-a q^{11}}^{-1} + \ul{\mr{1}_{-a q^{12}}^{-1}} \ .
}
Using the $q$-characters, we compute the poles of $\check{R}(z)$ and the corresponding kernels and cokernels. 
\begin{lem}
The poles of the $R$-matrix, the corresponding submodules and quotient modules are given by
\begin{center}\begin{tabular}{c@{\hspace{-1cm}} c@{\hspace{-1cm}} c}
Poles & Submodules & Quotient modules \\

$q^{2}$ & $\tl{L}_{\mr{1}_a\mr{1}_{aq^{-2}}} \cong L_{2\omega_1}\oplus L_{\omega_1}\oplus L_{\omega_0}$ & $\hspace{55pt}\tl{L}_{\mr{2}_{aq^{-1}}}\cong L_{\omega_2}\oplus L_{\omega_4}\oplus 2L_{\omega_1}\oplus L_{\omega_0}$ \\

$-q^{6}$ & $\hspace{29pt}\tl{L}_{\mr{1}_a\mr{1}_{-aq^{-6}}}\cong L_{2\omega_1}\oplus L_{\omega_2}\oplus 2L_{\omega_1}\oplus L_{\omega_0}$ & $\hspace{19pt} \tl{L}_{\mr{4}_{a^2q^{-6}}}\cong L_{\omega_4}\oplus L_{\omega_1}\oplus L_{\omega_0}$ \\

$q^{8}$ & $\hspace{63pt} \tl{L}_{\mr{1}_a\mr{1}_{aq^{-8}}}\cong L_{2\omega_1}\oplus L_{\omega_2}\oplus L_{\omega_4}\oplus 2L_{\omega_1}\oplus L_{\omega_0}$ & $\hspace{-9pt}\tl{L}_{\mr{1}_{-aq^{-4}}}\cong L_{\omega_1}\oplus L_{\omega_0}$ \\

$-q^{12}$ & $\hspace{56pt}\tl{L}_{\mr{1}_a\mr{1}_{-aq^{-12}}}\cong L_{2\omega_1}\oplus L_{\omega_2}\oplus L_{\omega_4}\oplus 3L_{\omega_1}\oplus L_{\omega_0}$ & $\hspace{-15pt} \tl{L}_{\scriptscriptstyle{1}}\cong L_{\omega_0}$\\
\end{tabular}\ .\end{center}
\p{
The proof is similar to the proof of Lemma \ref{poles Ato}.
}
\end{lem}

We choose a basis $\{v_i:1\le i\le 26\}\cup\{v_{27}\}$ for $L_{\omega_1}\oplus L_{\omega_0}$. A diagram of $L_{\omega_1}$ is given in \cite{DM25}. The vectors $v_{13}$, $v_{14}$ and $v_{27}$ are of weight zero. The action of $F_i$'s on the two dimensional $U_q\fk{sl}_2$-submodules inside $L_{\omega_1}$, involving $v_{12}, v_{13}, v_{14}, v_{15}$ with respect to $F_1$ and $v_{11}, v_{13}, v_{16}$ with respect to $F_2$ is given by 
$$F_1v_{12}=\frac{1}{\sqrt{[2]}}v_{13}+\frac{\sqrt{[3]}}{\sqrt{[2]}}v_{14}\ ,\> F_1v_{13}=\frac{1}{\sqrt{[2]}}v_{15}\ ,\> F_1v_{14}=\frac{\sqrt{[3]}}{\sqrt{[2]}}v_{15}\ ,\quad F_2v_{11}=\sqrt{[2]}v_{13}\ ,\> F_2v_{13}=\sqrt{[2]}v_{16}\ ,$$
All other $U_q\fk{sl}_2$-submodules in $L_{\omega_1}$ are one dimensional and $F_i$'s act by constant $1$. The $E_i$'s act in $L_{\omega_1}$ as transpose of $F_i$'s. The vector $v_{27}\in L_{\omega_0}\se\tl{L}_1(a)$, and is annihilated by all $F_i$'s and $E_i$'s.

A singular vectors of weight $2\omega_1$, respectively $\omega_2$, is chosen to be  $v_1\otimes v_1$, respectively $q\,v_1\otimes v_2-v_2\otimes v_1$.
A singular vector of weight $\omega_4$ is chosen to be 
$$q^3\,v_1\otimes v_7-q^2\,v_2\otimes v_6+q\,v_3\otimes v_4-q^{-1}\,v_4\otimes v_3+q^{-2}\,v_6\otimes v_2-q^{-3}\,v_7\otimes v_1\ .$$
We choose the three singular vectors of weight $\omega_1$ to be respectively
\bee{
& u_1=\frac{\sqrt{[3]}}{\sqrt{[4]}}\bigg(\frac{\sqrt{[2]}}{\sqrt{[3]}}q^6\,v_1\otimes v_{14}-q^{9/2}\,v_2\otimes v_{12} + q^{7/2}\,v_3\otimes v_{10}- q^{3/2}\,v_4\otimes v_{8} + q^{-1/2}\,v_5\otimes v_{6} \\
& + q^{1/2}\,v_{6}\otimes v_5  - q^{-3/2}\,v_{8}\otimes v_4 + q^{-7/2}\,v_{10}\otimes v_3 - q^{-9/2}\,v_{12}\otimes v_2+ \frac{\sqrt{[2]}}{\sqrt{[3]}}q^{-6}\,v_{14}\otimes v_1\bigg) \in L_{\omega_1}^{\otimes 2}\ ,\\ 
& u_2=v_1\otimes v_{27}\in L_{\omega_1}\otimes L_{\omega_0}\quad \text{and}\quad u_3=v_{27}\otimes v_1\in L_{\omega_0}\otimes L_{\omega_1}\ .
}
We choose the two singular vectors of weight $\omega_0$ to be respectively
$$w_1=\Big(\sum_{i=1}^{27} p_{i}^q\,v_i\otimes v_{27-i}\Big)+v_{13}\otimes v_{13}+v_{14}\otimes v_{14}\in L_{\omega_1}^{\otimes 2} \qquad \text{and} \qquad w_2=v_{27}\otimes v_{27}\in L_{\omega_0}^{\otimes 2}\ ,$$
where $p_i^q$ are as follows: $\{p_i^q: 1\le i\le 13\}$ is given by $\big\{q^{11}, -q^{10}, q^9, -q^7, q^5, q^6, -q^5, -q^4, q^3, q^2, -q, -q, 0\big\}$ and $p_i^q=p_{27-i}^{q^{-1}}$ for $14\le i\le 27$.

For $\la=2\omega_1,\omega_2,\omega_4,\omega_1,\omega_0$, let $P_\la^q$ be the projector onto the $U_q($F$_4)$-module $L_\la$ in the decomposition \eqref{tensorE6t}.

\begin{thm}
\label{thm:R E6t2}
In terms of projectors, we have 
\eq{\label{R proj E6t2}
\check{R}(z)=P_{2\omega_1}^q-q^{-2}\frac{1-q^{2}z}{1-q^{-2}z}\,P_{\omega_2}^q-q^{-8}\frac{(1-q^{2}z)(1+q^{6}z)}{(1-q^{-2}z)(1+q^{-6}z)}\,P_{\omega_4}^q+\frac{q^{-8}\,f_{\omega_1}(z)}{(1-q^{-2}z)(1+q^{-6}z)(1-q^{-8}z)}\otimes P_{\omega_1}^q \\
+\frac{q^{-14}\,f_{\omega_0}(z)}{(1-q^{-2}z)(1+q^{-6}z)(1-q^{-8}z)(1+q^{-12}z)}\otimes P_{\omega_0}^q\ ,
}
where the matrices $f_{\omega_1}(z)$ and $f_{\omega_0}(z)$ are given by 
$$f_{\omega_1}(z)=\begin{bmatrix}(q^{-3}-q^{3}z)(q^{-3}+\A\,z-q^{3}z^2) & \B\,z(1-z) & \B\,z(1-z)
\vspace{0.25cm}\\
\G\,z(1-z) & \B\,z(q^{6}-q^{-6}z) & (1-z)(q^{6}+\B\,z-q^{-6}z^2)
\vspace{0.25cm}\\
\G\,z(1-z) & (1-z)(q^{6}+\B\,z-q^{-6}z^2) & \B\,z(q^{6}-q^{-6}z)\end{bmatrix}\ ,$$
$$f_{\omega_0}(z)=\begin{bmatrix}q^{-12}+q^{-6}\,\zeta\,z+\xi\,z^2-q^{6}\,\zeta\,z^3+q^{12}z^4 & \eta\,z(1-z^2) \vspace{0.25cm}\\ \rho\,z(1-z^2) & q^{12}-q^{6}\,\zeta\,z+\xi\,z^2+q^{-6}\,\zeta\,z^3+q^{-12}z^4 \end{bmatrix}\ .$$
Here the constants $\A,\B,\G,\zeta,\xi,\eta,\rho\in\C(q)$ are given by
$$\A=\frac{[2]^{\mr{i}}\big([3]-[2]_{6}\big)}{[3]}\ ,\quad \B=\frac{[2]^{\mr{i}}\,[4]}{[3]}\ ,\quad \G=\frac{[2]^{\mr{i}}\,[2]\,[2]_6\,[7]}{[3]}\ ,\quad \zeta=\frac{[2]^{\mr{i}}\,[2]\,[7]}{[3]}\ ,$$
$$\eta=\frac{[2]^{\mr{i}}\,[4]}{[3]}\ ,\quad\rho=\frac{[2]^{\mr{i}}\,[3]_{3}^{\mr{i}}\,[2]_4\,[4]\,[13]}{[3]}\ ,\quad \xi=[2]_{14}-[2]_{12}-[2]_6+[2]_4-2\ .$$
\p{
In the expression of $\check{R}(z)$, the rational functions corresponding to the first three summands in \eqref{tensorE6t} are determined completely using the $q$-characters. Let $g_1(z)$ be the $3\times 3$ matrix and $g_2(z)$ be the $2\times 2$ matrix, corresponding to the last two summands $L_{\omega_1}$ and $L_{\omega_0}$ respectively.

The $3\times 3$ matrix $g_1(z)$ is determined (up to a sign) as follows. Using Lemma \ref{R0}, we get 
\eq{\label{e6g10}
g_1(0)=\begin{bmatrix} q^{-14} & 0 & 0\\ 0 & 0 & q^{-2} \\ 0 & q^{-2} & 0
\end{bmatrix}\ ,\quad g_1(\infty)=\begin{bmatrix} q^{14} & 0 & 0\\ 0 & 0 & q^{2} \\ 0 & q^{2} & 0
\end{bmatrix}\ .
}
From $q$-characters we know the poles of $g_1(z)$ and by Conjecture \ref{conj:simple poles} we presume that the poles are simple. Combining this and \eqref{e6g10} with $g_1(1)$ being zero on off-diagonal entries and that $g_1(z)$ commutes with the flip operator acting on singular vectors, see Lemma \ref{R and flip commute}, we get 
\bee{
g_1(z)=\frac{q^{-8}\,f_{\omega_1}(z)}{(1-q^{-2}z)(1+q^{-6}z)(1-q^{-8}z)}\ ,
}
where
\bee{
f_{\omega_1}(z)=\begin{bmatrix}
q^{-6}+\A_1 z+\A_2 z^2+q^{6}z^3 & \B z(1-z) & \B z(1-z) \\ 
\G z(1-z) & z(a_1+a_2z) & (1-z)(q^{6}+b z-q^{-6}z^2) \\
\G z(1-z) & (1-z)(q^{6}+b z-q^{-6} z^2) & z(a_1+a_2z)
\end{bmatrix}\ .
}
Since $g_1(1)$ is $1$ on the diagonal entries we have
\eq{\label{e6 g1 at 1}
a_1+a_2=[2]^{\mr{i}}\,[2]_3\,[2]_4^{\mr{i}}\ .
}
From $g_1(z)g_1(z^{-1})=\id$, we get
\eq{\label{e6 a1 a2}
a_1=-q^{12}a_2
}
and 
\eq{\label{e6 A1 A2}
a_2-\A_1-b=q^{-6}\ ,\quad a_1+b-\A_2=q^{6}\ .
}
The rank of $g_1(q^{-2})$ is $1$. This gives
\eq{\label{e6 g1 at q2}
q\,a_1+q^{-1}\,a_2=[2]^{\mr{i}}\,\big(b+[2]_8^{\mr{i}}\big)\ ,
}
and
\eq{\label{e6 g1 at q2 bg}
\big([2]^{\mr{i}}\big)^2\,\B\G=(q\,a_1+q^{-1}\,a_2)\,\big(q\,\A_1+q^{-1}\,\A_2+[2]_3\big)\ .
}
Now, using \eqref{e6 g1 at 1} and \eqref{e6 a1 a2} we get $a_1$ and $a_2$. Then \eqref{e6 g1 at q2} gives $b$. Then $\A_1$ and $\A_2$ are obtained using \eqref{e6 A1 A2}. Finally, the product $\B\G$ is obtained using \eqref{e6 g1 at q2 bg} and from the choice of singular vectors $u_1\in L_{\omega_1}^{\otimes 2}$, $u_2\in L_{\omega_1}\otimes L_{\omega_0}$ we have 
$$\frac{\G}{\B}=\frac{(u_1,u_1)}{(u_2,u_2)}=\frac{[2]\,[2]_6\,[7]}{[4]}\ .$$
This determines $f_{\omega_1}(z)$ up to the sign of $\B$ (or $\G$).

\medskip

The $2\times 2$ matrix $g_2(z)$ is determined (up to a sign) as follows. Using Lemma \ref{R0}, we get
\eq{\label{e6g20}
g_2(0)=\begin{bmatrix}
q^{-26} & 0\\ 0 & q^{-2}
\end{bmatrix}\ ,\quad g_2(\infty)=\begin{bmatrix}
q^{26} & 0\\ 0 & q^2
\end{bmatrix}\ .
}
From $q$-characters we know the poles of $g_2(z)$ and by Conjecture \ref{conj:simple poles} we presume that the poles are simple. Combining this and \eqref{e6g20} with $g_2(1)$ begin zero on off-diagonal entries we get
\bee{
g_2(z)=\frac{q^{-14}\,f_{\omega_0}(z)}{(1-q^{-2}z)(1+q^{-6}z)(1-q^{-8}z)(1+q^{-12}z)}\ ,
}
where
\bee{
f_{\omega_0}(z)=\begin{bmatrix}
q^{-12}+\zeta_1 z+\xi_1 z^2+\zeta_2 z^3+q^{12}z^4 & z(1-z)(\eta_1+\eta_2 z) \\ z(1-z)(\rho_1+\rho_2 z) & q^{12} +\zeta_3 z+\xi_2 z^2 + \zeta_4 z^3 +q^{-12} z^4
\end{bmatrix}\ .
}
Using $g_2(z)g_2(z^{-1})=\id$, we get
\bee{
\zeta_1=\zeta_4\ ,\quad \zeta_2=\zeta_3\ ,\quad \xi_1=\xi_2\ ,\quad \eta_1=\eta_2\ ,\quad \rho_1=\rho_2\ ,
}
so that
\bee{
f_{\omega_0}(z)=\begin{bmatrix}
q^{-12}+\zeta_1 z+\xi z^2+\zeta_2 z^3+q^{12}z^4 & \eta\, z(1-z^2) \\ \rho\, z(1-z^2) & q^{12} +\zeta_2 z+\xi z^2 + \zeta_1 z^3 +q^{-12} z^4
\end{bmatrix}\ .
}
Since $g_2(1)$ is $1$ on the diagonal entries we have
\eq{\label{e6 g2 at 1}
\zeta_1+\xi+\zeta_2+[2]_{12}=[2]^{\mr{i}}\,[2]_3\,[2]^{\mr{i}}_{4}\,[2]_{6}\ .
}
From $g_2(z)g_2(z^{-1})=\id$, now we get
\eq{\label{e6 g2 inv a}
q^{12}\,\zeta_1+q^{-12}\,\zeta_2=[2]^{\mr{i}}\,[2]\,[2]_3\,[2]_7^{\mr{i}}\ ,
}
\eq{\label{e6 g2 inv b}
q^{-12}\zeta_1+q^{12}\zeta_2+\xi(\zeta_1+\zeta_2)=-[2]^{\mr{i}}\,[2]\,[2]_3\,[2]_7^{\mr{i}}\,\big([2]_{14}-[2]_{6}+[2]_{4}-1\big)\ ,
}
\eq{\label{e6 g2 inv c}
\eta\rho=\zeta_1\zeta_2+\xi\,[2]_{12}+\big([2]_{20}-[2]_{18}+2[2]_{14}+[2]_{8}-2[2]_{6}+2[2]_{4}+[2]_{2}-4\big)\ .
}
Now, using \eqref{e6 g2 at 1}, \eqref{e6 g2 inv a} and \eqref{e6 g2 inv b} we get two solutions for each of $\zeta_1$, $\zeta_2$ and $\xi$. One of these solutions can be rejected  for the following reason.

If this solutions was the answer, then in the limit as $q$ goes to a primitive $24$-th root of unity, we observe that $L_{\omega_0}\otimes L_{\omega_0}=\C v_{27}\otimes v_{27}$ splits as a direct summand in $\tl{L}_1(z)\otimes\tl{L}_1$ for every $z$. That is possible only if $\tilde L_1(z)=L_{\omega_1}\oplus L_{\omega_0}$ as a $U_q\tl{\fk{g}}$-module for all $z$. It is easy to argue, see \eqref{E6t E0 action} below, that it is not the case when $q$ is a primitive root of unity of order $24$. We omit further details, and confirm our choice by checking that $R$-matrix does not commute with $E_0$ and $F_0$ for that choice of the solution.

After that we have a unique solution for $\zeta_1$, $\zeta_2$, $\xi$. Finally the product $\eta\,\rho$ is found using \eqref{e6 g2 inv c}, and from the choice of singular vectors $w_1\in L_{\omega_1}^{\otimes2}$, $w_2\in L_{\omega_0}^{\otimes2}$, we have 
$$\frac{\rho}{\eta}=\frac{(w_1,w_1)}{(w_2,w_2)}=[2]_4\,[3]_3^{\mr{i}}\,[13]\ .$$
This determines $f_{\omega_0}(z)$ up to the sign of $\eta$ (or $\rho$).

To fix the signs of $\B$ in $f_{\omega_1}(z)$ and $\eta$ in $f_{\omega_0}(z)$, we use the $E_0$ action. Namely, to determine the sign of $\B$ we apply both sides of the commutation relation in \eqref{R and E0 commutation relation} to $v_1\otimes v_1$ and compare the coefficients of $v_1\otimes v_{27}$ on the two sides. To determine the sign of $\eta$ we apply both sides of \eqref{R and E0 commutation relation} to $v_1\otimes v_{27}$ and compare coefficients of $v_{27}\otimes v_{27}$ on the two sides.

One directly checks that the $R$-matrix commutes with the action of $E_0$ and $F_0$, where  
\bee{
K_0 =\ & q^{-2}\,E_{11}\ +\ q^{-1}\,\sum_{i=2}^5\big(E_{ii}+E_{16-2i,16-2i}\big)\ +\ \sum_{i=3}^5\big(E_{2i+1,2i+1}+E_{26-2i,26-2i}\big)\\ 
& +\ E_{13,13}\ +\ E_{27,27}\ +\ q\,\sum_{i=2}^5\big(E_{11+2i,11+2i} + E_{27-i,27-i}\big)\ +\ q^2\,E_{26,26}\ ,
}
\eq{\label{E6t E0 action}
E_0(a)=\ & a\frac{\sqrt{[2]}}{\sqrt{[3]}}\,\big(E_{14,1}+E_{26,14}\big)\ + a\frac{\sqrt{[4]}}{\sqrt{[3]}}\,\big(E_{27,1}+E_{26,27}\big)\ + a\,\sum_{i=2}^5\big(E_{11+2i,i}+E_{27-i,16-2i}\big)\ ,
}
and $F_0(a)$ is the transpose of $a^{-2}E_0(a)$.
}\end{thm}

In the rational case, we recover the untwisted type E$_{6}^{(1)}$ rational $R$-matrix in Corollary 5.3 in \cite{DM25} as follows. Let $\check{R}(u)$ be the rational $R$-matrix obtained after substituting $z=q^{2u}$ in \eqref{R proj E6t2} and taking the $q\to 1$ limit. Let $T:\C^{27}\to\C^{27}$ be a linear map given by $T(v_i)=v_i$, for $1\le i\le 12,\ i\ne 7,8$, $T(v_7)=v_8$, $T(v_8)=v_7$, $T(v_i)=v_{i+1}$ for $15\le i\le 26,\ i\ne 19,20$, $T(v_{19})=v_{21}$, $T(v_{20})=v_{20}$, $T(v_{13})=\frac{1}{\sqrt{2}}(v_{13}+v_{14})$, $T(v_{14})=\frac{-1}{\sqrt{6}}(v_{13} - v_{14} - 2\,v_{15})$, and $T(v_{27})=\frac{-1}{\sqrt{3}}(v_{13}-v_{14}+v_{15})$. Then the matrix $(T\otimes T)\check{R}(u)(T\otimes T)^{-1}$ is the untwisted type E$_{6}^{(1)}$ rational $R$-matrix.



\subsection{Type D\texorpdfstring{$_4^{(3)}$}{2}}

The $8$-dimensional $U_q\big($D$_4^{(3)}\big)$-module $\tl{L}_1(a)$ restricted to $U_q($G$_2)$ is isomorphic to $L_{\omega_1}\oplus L_{\omega_0}$. \\
As $U_q($G$_2)$-modules, we have 
\eq{\label{tensorDt3}
\big(\underbracket[0.1ex]{L_{\omega_1}}_{7}\oplus \underbracket[0.1ex]{L_{\omega_0}}_{1}\big)^{\otimes 2}\cong \underbracket[0.1ex]{L_{2\omega_1}}_{27}\oplus \underbracket[0.1ex]{L_{\omega_2}}_{14}\oplus \,3\underbracket[0.1ex]{L_{\omega_1}}_{7}\oplus \, 2\underbracket[0.1ex]{L_{\omega_0}}_{1}\ .
}
In the $q\to1$ limit, $L_{2\omega_1}\oplus L_{\omega_0}\mapsto \cl{S}^2(L_{\omega_1})$ and $L_{\omega_2}\oplus L_{\omega_1}\mapsto \Lambda^2(L_{\omega_1})$.

The $q$-character of $\tl{L}_{\mr{1}_a}$ has $8$ terms and there are $2$ weight zero terms (shown in box):
\bee{
\chi_q(\mr{1}_a)=\mr{1}_a + \ul{\mr{1}_{aq^2}^{-1} \mr{2}_{aq}} + \mr{1}_{\mr{j}aq^2} \mr{1}_{\mr{j}^2aq^2} \mr{2}_{aq^3}^{-1} + \boxed{\ul{\mr{1}_{\mr{j}^2aq^4}^{-1} \mr{1}_{\mr{j}aq^2}}} + \boxed{\ul{\mr{1}_{\mr{j}aq^4}^{-1} \mr{1}_{\mr{j}^2aq^2}}} + \mr{1}_{\mr{j}aq^4}^{-1} \mr{1}_{\mr{j}^2aq^4}^{-1} \mr{2}_{aq^3} + \mr{1}_{aq^4} \mr{2}_{aq^5}^{-1} + \ul{\mr{1}_{aq^6}^{-1}}\ .
}
Using the $q$-characters, we compute the poles of $\check{R}(z)$ and the corresponding kernels and cokernels.

\begin{lem}
The poles of the $R$-matrix, the corresponding submodules and quotient modules are given by
\begin{center}\begin{tabular}{c c c }
Poles & Submodules & Quotient modules \\

$q^{2}$ & $\tl{L}_{\mr{1}_a\mr{1}_{aq^{-2}}}\cong L_{2\omega_1}\oplus L_{\omega_1}\oplus L_{\omega_0}$ & $\hspace{36pt} \tl{L}_{\mr{2}_{a^3q^{-3}}}\cong L_{\omega_2}\oplus 2L_{\omega_1}\oplus L_{\omega_0}$ \\

$\mr{j}q^{4}$ & $\hspace{31pt} \tl{L}_{\mr{1}_a\mr{1}_{\mr{j}aq^{-4}}}\cong L_{2\omega_1}\oplus L_{\omega_2}\oplus 2L_{\omega_1}\oplus L_{\omega_0}$ & $\tl{L}_{\mr{1}_{\mr{j}^2aq^{-2}}}\cong L_{\omega_1}\oplus L_{\omega_0}$ \\

$\mr{j}^2q^{4}$ & $\hspace{28pt} \tl{L}_{\mr{1}_a\mr{1}_{\mr{j}^2aq^{-4}}}\cong L_{2\omega_1}\oplus L_{\omega_2}\oplus 2L_{\omega_1}\oplus L_{\omega_0}$ & $\hspace{4pt} \tl{L}_{\mr{1}_{\mr{j}aq^{-2}}}\cong L_{\omega_1}\oplus L_{\omega_0}$ \\

$q^{6}$ & $\hspace{34pt} \tl{L}_{\mr{1}_a\mr{1}_{aq^{-6}}}\cong L_{2\omega_1}\oplus L_{\omega_2}\oplus 3L_{\omega_1}\oplus L_{\omega_0}$ & $\hspace{-5pt} \tl{L}_{\scriptscriptstyle{1}}\cong L_{\omega_0}$\\
\end{tabular}\end{center}
\p{
The proof is similar to the proof of Lemma \ref{poles Ato}.
}
\end{lem}

We choose a basis $\{v_i:1\le i\le 7\}\cup\{v_{8}\}$ for $L_{\omega_1}\oplus L_{\omega_0}$. The vectors $v_4$ and $v_8$ are of weight zero. The action of $F_1, F_2$ in $L_{\omega_1}$ is given by
$$F_1v_i=v_{i+1}\ ,i=1,6\ ,\quad F_2v_i=v_{i+1}\ ,i=2,5\ ,\quad F_1v_i=\sqrt{[2]}\,v_{i+1}\ ,i=3,4\ .$$
The generators $E_1, E_2$ act in $L_{\omega_1}$ as transpose of $F_1, F_2$ respectively. The vector $v_{8}\in L_{\omega_0}\se\tl{L}_1(a)$, and is annihilated by $F_1, F_2$ and $E_1, E_2$.

A singular vectors of weight $2\omega_1$, respectively $\omega_2$, is chosen to be  $v_1\otimes v_1$, respectively $q\,v_1\otimes v_2-v_2\otimes v_1$.
The three singular vectors of weight $\omega_1$ are chosen to be respectively
\bee{
& u_1=\frac{1}{\sqrt{[3]}}\bigg(q^3\,v_1\otimes v_{4}-q^{3/2}\,\sqrt{[2]}\,v_2\otimes v_{3} + q^{-3/2}\,\sqrt{[2]}\,v_3\otimes v_{2} -q^{-3}\,v_{4}\otimes v_1\bigg) \in L_{\omega_1}^{\otimes 2}\ ,\\ 
& u_2=v_1\otimes v_{8}\in L_{\omega_1}\otimes L_{\omega_0}\quad \text{and}\quad u_3=v_{8}\otimes v_1\in L_{\omega_0}\otimes L_{\omega_1}\ .
}
The two singular vectors of weight $\omega_0$ are chosen to be respectively
$$w_1=\sum_{i=1}^7\,p_i^q\,v_i\otimes v_{8-i}\in L_{\omega_1}^{\otimes 2} \qquad \text{and} \qquad w_2=v_{8}\otimes v_{8}\in L_{\omega_0}^{\otimes 2}\ ,$$
where $p_i^q$ are given by $\{q^5, -q^4, q, -1, q^{-1}, -q^{-4}, q^{-5}\}$.

For $\la=2\omega_1,\omega_2,\omega_1,\omega_0$, let $P_\la^q$ be the projector onto the $U_q($G$_2)$-module $L_\la$ in the decomposition \eqref{tensorDt3}.

\begin{thm}
\label{thm:R D4t3}
In terms of projectors, we have
\eq{\label{R proj D4t3}
\check{R}(z)=P_{2\omega_1}^q-q^{-2}\frac{1-q^{2}z}{1-q^{-2}z}\,P_{\omega_2}^q + \frac{q^{-5}\,f_{\omega_1}(z)}{(1-q^{-2}z)(1+q^{-4}z+q^{-8}z^2)}\otimes P_{\omega_1}^q \\
+\frac{q^{-8}\,f_{\omega_0}(z)}{(1-q^{-2}z)(1-q^{-6}z)(1+q^{-4}z+q^{-8}z^2)}\otimes P_{\omega_0}^q\ ,
}
where the matrices $f_{\omega_1}(z)$ and $f_{\omega_0}(z)$ are given by
$$f_{\omega_1}(z)=\begin{bmatrix}-q^{-3}-q^{-2}\,\A\,z+q^{2}\,\A\,z^2+q^{3}z^3 & \B\,z(1-z) & \B\,z(1-z)\vspace{0.25cm}\\ \G\,z(1-z) & \B\,z(q^{3}+q^{-3}z) & (1-z)(q^{3}+\kappa\,z+q^{-3}z^2) \vspace{0.25cm}\\ \G\,z(1-z) & (1-z)(q^{3}+\kappa\,z+q^{-3}z^2) & \B\,z(q^{3}+q^{-3}z)\end{bmatrix}\ ,$$
$$f_{\omega_0}(z)=\begin{bmatrix}q^{-6}-q^{-3}\,\zeta\,z+\xi\,z^2-q^{3}\,\zeta\,z^3+q^{6}z^4 & \eta\,z(1-z^2) \vspace{0.25cm}\\ \rho\,z(1-z^2) & q^{6}-q^{3}\,\zeta\,z+\xi\,z^2-q^{-3}\,\zeta\,z^3+q^{-6}z^4 \end{bmatrix}\ .$$
Here the constants $\A,\B,\G,\kappa,\zeta,\xi,\eta,\rho\in\C(q)$ are given by
$$\A=\frac{[2]^{\mr{i}}\,[2]_{3}^{\mr{i}}}{[2]}\ ,\> \B=\frac{[2]_{3}^{\mr{i}}}{[2]}\ ,\> \G=\frac{[2]_{3}^{\mr{i}}\,[2]_4}{[2]}\ ,\> \kappa=\frac{[2]_{2}}{[2]}\ ,\> \zeta=\frac{[2]_{4}}{[2]}\ ,\>\eta=\frac{[2]_{3}^{\mr{i}}}{[2]}\ ,\> \rho=\frac{[2]_{3}^{\mr{i}}\,[3]_{2}^{\mr{i}}\,[7]}{[2]}\ ,\> \xi=[2]^{\mr{i}}\,[2]_3^{\mr{i}}\,[2]_4\ .$$
\p{
In the expression of $\check{R}(z)$, the rational functions corresponding to the first three summands in \eqref{tensorDt3} are determined completely using the $q$-characters. Let $g_1(z)$ be the $3\times 3$ matrix and $g_2(z)$ be the $2\times 2$ matrix, corresponding to the last two summands $L_{\omega_1}$ and $L_{\omega_0}$ respectively.

The $3\times 3$ matrix $g_1(z)$ is determined (up to a sign) as follows. Using Lemma \ref{R0}, we get 
\eq{\label{d4g10}
g_1(0)=\begin{bmatrix} -q^{-8} & 0 & 0\\ 0 & 0 & q^{-2} \\ 0 & q^{-2} & 0
\end{bmatrix}\ ,\quad g_1(\infty)=\begin{bmatrix} -q^{8} & 0 & 0\\ 0 & 0 & q^{2} \\ 0 & q^{2} & 0
\end{bmatrix}\ .
}
From $q$-characters we know the poles of $g_1(z)$ and by Conjecture \ref{conj:simple poles} we presume that the poles are simple. Combining this and \eqref{d4g10} with $g_1(1)$ being zero on off-diagonal entries and that $g_1(z)$ commutes with the flip operator acting on singular vectors, see Lemma \ref{R and flip commute}, we get
\bee{
g_1(z)=\frac{q^{-5}\,f_{\omega_1}(z)}{(1-q^{-2}z)(1+q^{-4}z+q^{-8}z^2)}\ ,
}
where 
\bee{
f_{\omega_1}(z)=\begin{bmatrix}
-q^{-3}+\A_1 z+\A_2 z^2+q^{3}z^3 & \B z(1-z) & \B z(1-z) \\ 
\G z(1-z) & z(a_1+a_2z) & (1-z)(q^{3}+b z+q^{-3}z^2) \\
\G z(1-z) & (1-z)(q^{3}+b z+q^{-3} z^2) & z(a_1+a_2z)
\end{bmatrix}\ .
}
Since $g_1(1)$ is $1$ on the diagonal entries we have
\eq{\label{d4 g1 at 1}
a_1+a_2=[2]^{\mr{i}}\,[3]_2\ .
}
From $g_1(z)g_1(z^{-1})=\id$, we get
\eq{\label{d4 a1 a2}
a_1=q^{6}a_2
}
and 
\eq{\label{d4 A1 A2}
\A_1-a_2+b=q^{-3}\ ,\quad \A_2-a_1-b=-q^{3}\ .
}
The rank of $g_1(q^{-2})$ is $1$. This gives
\eq{\label{d4 g1 at q2}
q\,a_1+q^{-1}\,a_2=[2]^{\mr{i}}\,\big(b+[2]_5\big)\ ,
}
and
\eq{\label{d4 g1 at q2 bg}
\big([2]^{\mr{i}}\big)^2\,\B\,\G=(q\,a_1+q^{-1}\,a_2)(q\,\A_1+q^{-1}\,\A_2)\ .
}
Now, using \eqref{d4 g1 at 1} and \eqref{d4 a1 a2} we get $a_1$ and $a_2$. Then \eqref{d4 g1 at q2} gives $b$. Then $\A_1$ and $\A_2$ are obtained using \eqref{d4 A1 A2}. Finally, the product $\B\G$ is obtained using \eqref{d4 g1 at q2 bg}. From the choice of singular vectors $u_1\in L_{\omega_1}^{\otimes2}$ and $u_2\in L_{\omega_1}\otimes L_{\omega_0}$, we have
$$\frac{\G}{\B}=\frac{(u_1,u_1)}{(u_2,u_2)}=[2]_4\ .$$
This determines $f_{\omega_1}(z)$ up to the sign of $\B$ (or $\G$).

\medskip

The $2\times 2$ matrix $g_2(z)$ is determined (up to a sign) as follows. Using Lemma \ref{R0}, we get
\eq{\label{d4g20}
g_2(0)=\begin{bmatrix}
q^{-14} & 0\\ 0 & q^{-2}
\end{bmatrix}\ ,\quad g_2(\infty)=\begin{bmatrix}
q^{14} & 0\\ 0 & q^2
\end{bmatrix}\ .
}
From $q$-characters we know the poles of $g_2(z)$ and by Conjecture \ref{conj:simple poles} we presume that the poles are simple. Combining this and \eqref{d4g20} with $g_2(1)$ begin zero on off-diagonal entries we get
\bee{
g_2(z)=\frac{q^{-8}\,f_{\omega_0}(z)}{(1-q^{-2}z)(1-q^{-6}z)(1+q^{-4}z+q^{-8}z^2)}\ ,
}
where 
\bee{
f_{\omega_0}(z)=\begin{bmatrix}
q^{-6}+\zeta_1 z+\xi_1 z^2+\zeta_2 z^3+q^{6}z^4 & z(1-z)(\eta_1+\eta_2 z) \\ z(1-z)(\rho_1+\rho_2 z) & q^{6} +\zeta_3 z+\xi_2 z^2 + \zeta_4 z^3 +q^{-6} z^4
\end{bmatrix}\ .
}
Using $g_2(z)g_2(z^{-1})=\id$, we get
\bee{
\zeta_1=\zeta_4\ ,\quad \zeta_2=\zeta_3\ ,\quad \xi_1=\xi_2\ ,\quad \eta_1=\eta_2\ ,\quad \rho_1=\rho_2\ ,
}
so that
\bee{
f_{\omega_0}(z)=\begin{bmatrix}
q^{-6}+\zeta_1 z+\xi z^2+\zeta_2 z^3+q^{6}z^4 & \eta\, z(1-z^2) \\ \rho\, z(1-z^2) & q^{6} +\zeta_2 z+\xi z^2 + \zeta_1 z^3 +q^{-6} z^4
\end{bmatrix}\ .
}
Since $g_2(1)$ is $1$ on the diagonal entries we have
\eq{\label{d4 g2 at 1}
\zeta_1+\xi+\zeta_2+[2]_{6}=[2]^{\mr{i}}\,[2]^{\mr{i}}_3\,[3]_{2}\ .
}
From $g_2(z)g_2(z^{-1})=\id$, now we get
\eq{\label{d4 g2 inv a}
q^{6}\,\zeta_1+q^{-6}\,\zeta_2=-[2]_4\,[3]^{\mr{i}}\ ,
}
\eq{\label{d4 g2 inv b}
q^{-6}\,\zeta_1 + q^{6}\,\zeta_2 + \xi\,(\zeta_1+\zeta_2) = -[2]_4\,[3]^{\mr{i}}\,\big([2]_{8}-[2]_{2}+1\big)\ ,
}
\eq{\label{d4 g2 inv c}
\eta\,\rho=\zeta_1\zeta_2+\xi\,[2]_{6}+[2]_{10}-2\,[2]_{8}+[2]_{6}-[2]_{4}+2\,[2]_{2}-3\ .
}
Now, using \eqref{d4 g2 at 1}, \eqref{d4 g2 inv a} and \eqref{d4 g2 inv b} we get two solutions for each of $\zeta_1$, $\zeta_2$ and $\xi$, out of which one is rejected because the $q\to 1$ limit does not exist in that case. After that we have a unique solution for $\zeta_1$, $\zeta_2$, $\xi$. Finally, the product $\eta\,\rho$ is found using \eqref{d4 g2 inv c}. From the choice of singular vectors $w_1\in L_{\omega_1}^{\otimes 2}$ and $w_2\in L_{\omega_0}^{\otimes2}$, we have
$$\frac{\rho}{\eta}=\frac{(w_1,w_1)}{(w_2,w_2)}=[3]_2^{\mr{i}}\,[7]\ .$$
This determines $f_{\omega_0}(z)$ up to the sign of $\eta$ (or $\rho$).

\medskip

To fix the signs of $\B$ in $f_{\omega_1}(z)$ and $\eta$ in $f_{\omega_0}(z)$, we use the $E_0$ action. Namely, to determine the sign of $\B$ we apply both sides of the commutation relation in \eqref{R and E0 commutation relation} to $v_1\otimes v_1$ and compare the coefficients of $v_1\otimes v_8$ on the two sides. To determine the sign of $\eta$ we apply both sides of \eqref{R and E0 commutation relation} to $v_1\otimes v_8$ and compare coefficients of $v_8\otimes v_8$ on the two sides.

One can directly check that the $R$-matrix commutes with the action of $E_0$ and $F_0$, where  
\bee{
K_0 = q^{-2}\,E_{11}\ +\ q^{-1}\,\big(E_{22}+E_{33}\big)\ +\ \big(E_{44}+E_{88}\big)\ +\ q\,\big(E_{55} + E_{66}\big)\ +\ q^2\,E_{77}\ ,
}
\bee{
E_0(a)=\ & a\frac{1}{\sqrt{[2]}}\,\big(E_{41}+E_{74}\big)\ + a\frac{\sqrt{[3]}}{\sqrt{[2]}}\,\big(E_{81}+E_{78}\big)\ + a\,\big(E_{52}+E_{63}\big)\ ,
}
and $F_0(a)$ is the transpose of $a^{-2}E_0(a)$.
}\end{thm}

In the rational case, we recover the untwisted type D$_{4}^{(1)}$ rational $R$-matrix in Corollary 4.13 in \cite{DM25} as follows. Let $\check{R}(u)$ be the rational $R$-matrix obtained after substituting $z=q^{2u}$ in \eqref{R proj D4t3} and taking the $q\to 1$ limit. Let $T:\C^8\to\C^8$ be a linear map given by $T(v_i)=v_i$, for $1\le i\le 3$, $T(v_i)=v_{i+1}$ for $5\le i\le 7$, $T(v_{4})=v_{4}+\frac{1}{2}v_{5}$, and $T(v_{8})=\mr{i}\,v_{4}-\frac{\mr{i}}{2}v_{5}$, where $\mr{i}$ is the primitive second root of unity. Then the matrix $(T\otimes T)\check{R}(u)(T\otimes T)^{-1}$ is the untwisted type D$_{4}^{(1)}$ rational $R$-matrix.

\bigskip

{\bf Acknowledgments.\ }
The authors are partially supported by Simons Foundation grant number \#709444.


\bigskip

\begin{thebibliography}{0000000}
\bibitem[B85]{B85} V. Bazhanov, \textit{Trigonometric solutions of triangle equations and classical Lie algebras}, Physics Letters B \tb{159} (1985), no. 4, 321-324.
\bibitem[C00]{C00} V. Chari, \textit{Braid group actions and tensor products}, International Mathematics Research Notices (2000),  no. 7, 357–382.
\bibitem[CP98]{CP98} V. Chari and A. Pressley, \ti{Twisted quantum affine algebras}, Comm. Math. Phys. \tb{196} (1998), no. 2, 461-476.
\bibitem[DM25]{DM25} K. Dahiya and E. Mukhin, {\it Intertwiners of representations of untwisted quantum affine algebras and Yangians revisited}, arXiv:2503.09845 [math.QA].
\bibitem[Da98]{Da98} I. Damiani, \textit{The R-matrix for the (twisted) quantum affine algebras}, Representations and quantizations, Shanghai, (1998), 89-144, China High. Educ. Press, Beijing, 2000. 
\bibitem[Da14]{Da14} I. Damiani, \textit{Drinfeld Realization of Affine Quantum Algebras: the Relations}, Publ. Res. Inst. Math. Sci. \tb{48} (2012), no. 3, 661-733.
\bibitem[DGZ94]{DGZ94} G. Delius, M. Gould, and Y.-Z. Zhang, \textit{On the construction of trigonometric solutions of the Yang–Baxter equation}, Nuclear Physics B \tb{432} (1994), no. 1, 377-403. 
\bibitem[DGZ96]{DGZ96} G. Delius, M. Gould, and Y.-Z. Zhang, \textit{Twisted quantum affine algebras and solutions to the Yang-Baxter equation}, Int. J. Mod. Phys. A \textbf{11} (1996), 3415-3438. 
\bibitem[Dr87]{Dr87} V. Drinfeld, \textit{A new realization of Yangians and of quantum affine algebras}, Sov. Math. Dokl. \tb{36} (1987), 212-216.
\bibitem[FM01]{FM1} E. Frenkel and E. Mukhin, \textit{Combinatorics of q-characters of finite-dimensional representations of quantum affine algebras}, Commun. Math. Phys. \tb{216} (2001), 23-57.
\bibitem[GMW96]{GMW96} G. Gandenberger, N. MacKay, and G. Watts, \ti{Twisted algebra R-matrices and S-matrices for B$_n^{(1)}$ affine Toda solitons and their bound states}, Nucl.Phys. B \tb{465} (1996), 329-349.
\bibitem[H10]{H10} D. Hernandez, \textit{Kirillov–Reshetikhin Conjecture: The General Case}, International Mathematics Research Notices \tb{2010},  (2010), no. 1,  149-193. 
\bibitem[J86]{J86} M. Jimbo, \textit{Quantum R matrix for the generalized Toda system}, Comm. Math. Phys. \tb{102} (1986), no. 4, 537-547. 
\bibitem[(KMN)$^2$92]{KKMMNN92} S. Kang, M. Kashiwara, K. Misra, T. Miwa, T. Nakashima, A. Nakayashiki, \ti{Perfect crystals of quantum affine Lie algebras}, Duke Mathematical Journal, \tb{68} (1992), no. 3, 499-607.
\bibitem[KMOY06]{KMOY06} M. Kashiwara, K. Misra, M. Okado, D. Yamada, \ti{Perfect crystals for $U_q(D_4^{(3)})$}, Journal of Algebra \tb{317} (2006), 392-423.
\bibitem[KR90]{KR90} A. Kirillov, N. Reshetikhin, \ti{$q$-Weyl group and a multiplicative formula for universal $R$-matrices}, Comm. Math. Phys. \tb{134} (1990), 421-31.
\bibitem[K90]{K90} A. Kuniba, \textit{Quantum R matrix for G$_2$ and a solvable 175-vertex model}, Journal of Physics A: Mathematical and General \tb{23} (1990), no. 8, 1349-1362. 
\bibitem[LS90]{LS90} S. Levendorskii, Ya. Soibelman, \ti{Some applications of quantum Weyl groups}, J. Geom. Phys. \tb{7} (1990), 241-54.
\bibitem[M90]{M90} Z.-Q. Ma,  {\it The spectrum-dependent solutions to the Yang-Baxter equation for quantum E$_6$ and E$_7$}, J. Phys. A {\bf 23} (1990), no. 23, 5513-5522.
\bibitem[M91]{M91} Z.-Q. Ma, \textit{The embedding $e_0$ and the spectrum-dependent R-matrix for q-F$_4$}, J. Phys. A: Math. Gen. {\bf 24} (1991), 433-449.
\bibitem[MY14]{MY14} E. Mukhin and C. Young, \textit{Affinization of category O for quantum groups}, Transactions of the American Mathematical Society \tb{366} (2014), no. 9, 4815-4847.
\bibitem[O86]{O86} E. I. Ogievetsky, \textit{Factorised $S$-matrix with $G_2$ symmetry}, J. Phys. G: Nucl. Phys. \tb{12} (1986), L105-L108.
\bibitem[ZJ20]{ZJ20} P. Zinn-Justin, \textit{The trigonometric E$_8$ $R$-matrix},  Lett. Math. Phys. \tb{110} (2020), no. 12, 3279-3305.
\end{thebibliography}
\end{document}